\documentclass[a4paper]{article}%
\usepackage{amsmath}
\usepackage{amsfonts}
\usepackage{amssymb}
\usepackage{graphicx}
\usepackage{hyperref}%
\providecommand{\U}[1]{\protect\rule{.1in}{.1in}}
\newtheorem{theorem}{Theorem}[section]

\newtheorem{corollary}[theorem]{Corollary}

\newtheorem{lemma}[theorem]{Lemma}

\newtheorem{proposition}[theorem]{Proposition}
\newtheorem{remark}[theorem]{Remark}

\numberwithin{equation}{section}
\newenvironment{proof}[1][Proof]{\noindent\textbf{#1.} }{\ \rule{0.5em}{0.5em}}
\begin{document}

\title{A convective $1$-homogeneous $p$-Laplacian problem}
\author{Grey Ercole\\{\small Universidade Federal de Minas Gerais} \\{\small Belo Horizonte, MG, 30.123-970, Brazil}\\{\small grey@mat.ufmg.br}}
\maketitle

\begin{abstract}
Let $\Omega\subset\mathbb{R}^{N}$ be a bounded, smooth domain, and let
$a:\overline{\Omega}\rightarrow\mathbb{R}$ be a continuous function that may
change sign. We study the $1$-homogeneous problem
\[
\left\{
\begin{array}
[c]{lll}
-\operatorname{div}\left(  \left\vert \nabla u\right\vert ^{p-2}\nabla
u\right)  =\lambda\left\vert u\right\vert ^{p-2}u+a(x)\left\vert \nabla
u\right\vert ^{p-1} & \mathrm{in} & \Omega\\
u=0 & \mathrm{on} & \partial\Omega,
\end{array}
\right.
\]
where $1<p<\infty.$ We prove the existence of $\lambda_{a}>0$ for which this
problem,  with $\lambda = \lambda_a$, has a positive weak solution $u_{a}\in W_{0}^{1,p}(\Omega)\cap
C^{1}(\overline{\Omega}).$ The significance of this result lies in the fact
that no restriction is imposed on the size of $\left\Vert a\right\Vert
_{\infty}.$ We also obtain lower and upper bounds for $\lambda_{a}$ in terms of the positive and negative parts of $a$ as well as the circumradius and
inradius of $\Omega.$

\end{abstract}



\section{Introduction}

In this paper we study the existence of positive weak solutions to the
following nonlinear eigenvalue-type problem involving a convective term%
\begin{equation}
\left\{
\begin{array}
[c]{lll}%
-\Delta_{p}u=\lambda\left\vert u\right\vert ^{p-2}u+a(x)\left\vert \nabla
u\right\vert ^{p-1} & \mathrm{in} & \Omega\\
u=0 & \mathrm{on} & \partial\Omega,
\end{array}
\right.  \label{Pw}%
\end{equation}
where $\Omega\subset\mathbb{R}^{N}$ ($N\geq2$) is a bounded domain with a
$C^{2}$ boundary, $\Delta_{p}u:=\operatorname{div}(\left\vert \nabla
u\right\vert ^{p-2}\nabla u)$ is the $p$-Laplacian operator, $\lambda>0,$ and
$a\in C(\overline{\Omega})$ is a drift coefficient with no restriction on its sign.

We observe that (\ref{Pw}) is $1$-homogeneous: if $u$ is a nontrivial weak
solution and $k\not =0$ is a constant, then $ku$ is also a nontrivial weak
solution. Furthermore, this problem does not have a variational structure
because of the convective term.

The literature on $1$-homogeneous Dirichlet problems involving the
$p$-Laplacian and a sign-changing convective term is rather limited.

In \cite{ACG}, Anane, Chakrone and Gossez studied the problem%
\begin{equation}
\left\{
\begin{array}
[c]{lll}%
-\Delta_{p}u=\lambda\left\vert u\right\vert ^{p-2}u+(z.\nabla u)\left\vert
\nabla u\right\vert ^{p-2} & \mathrm{in} & \Omega\\
u=0 & \mathrm{on} & \partial\Omega,
\end{array}
\right.  \label{Pdot}%
\end{equation}
where $\lambda\in\mathbb{R}$ and $z\in\mathbb{R}^{N}$ (the dot
\texttt{"$\cdot$"} stands for the scalar product in $\mathbb{R}^{N}$). This
problem is equivalent to the eigenvalue problem%
\begin{equation}
\left\{
\begin{array}
[c]{lll}%
-\operatorname{div}(e^{z\cdot x}\left\vert \nabla u\right\vert ^{p-2}\nabla
u)=\lambda e^{z\cdot x}\left\vert u\right\vert ^{p-2}u & \mathrm{in} &
\Omega\\
u=0 & \mathrm{on} & \partial\Omega,
\end{array}
\right.  \label{dot}%
\end{equation}
which is variational. This fact allowed the authors to apply the
Ljusternik-Schnirelmann minimax method and show that (\ref{Pdot}) enjoys
essentially the same spectral properties as the Dirichlet $p$-Laplacian (i.e. $z=0$). In
particular, they proved that for each $z\in\mathbb{R}^{N}$ there exist
$\lambda_{z}>0$ and a positive function $u_{z}\in W_{0}^{1,p}(\Omega)$ such
that $u_{z}$ is a weak solution to (\ref{Pdot}) with $\lambda=\lambda_{z}.$
Moreover, $\lambda_{z}$ is simple: if $u\in W_{0}^{1,p}(\Omega)$ is another
nontrivial weak solution to (\ref{Pdot}) with $\lambda=\lambda_{z},$ then
$u=ku_{z}$ for some constant $k\not =0.$

In \cite{GST}, Garc\'{\i}a-Mel\'{\i}an, Sabina de Lis and Tak\'{a}\v{c}
studied the problem%
\begin{equation}
\left\{
\begin{array}
[c]{lll}%
-\Delta_{p}u+B(x,\nabla u)=\lambda\left\vert u\right\vert ^{p-2}u &
\mathrm{in} & \Omega\\
u=0 & \mathrm{on} & \partial\Omega,
\end{array}
\right.  \label{B}%
\end{equation}
where the continuous function $B:\overline{\Omega}\times\mathbb{R}%
^{N}\rightarrow\mathbb{R}$ is assumed to be $(p-1)$-homogeneous and odd with
respect to the second variable: $B(x,s\xi)=s\left\vert s\right\vert
^{p-2}B(x,\xi)$ for a.e. $x\in\Omega$ and for all $\xi\in\mathbb{R}^{N}$ and
$s\in\mathbb{R}$. A prototype of $B,$ which includes
(\ref{Pdot}) as a particular case of (\ref{B}), is $B(x,\xi)=(A(x)\cdot
\xi)\left\vert \xi\right\vert ^{p-2},$ with $A\in C(\Omega)^{N}.$ Garc\'{\i}a-Mel\'{\i}an, Sabina de Lis and Tak\'{a}\v{c}
combined variational arguments with comparison techniques and proved that
(\ref{B}) has a unique eigenvalue $\lambda_{1}>0$ associated with a positive
eigenfunction $u_{1}\in W_{0}^{1,p}(\Omega).$ They also proved that
$\lambda_{1}$ is simple and isolated. We remark that (\ref{B}) does not
include (\ref{Pw}) since the function $(x,\xi)\mapsto-a(x)\left\vert
\xi\right\vert ^{p-1}$ is even with respect to $\xi.$

The unidimensional version of \eqref{B}, in the particular case where $B(x,\xi)=c\left\vert \xi\right\vert ^{p-2}\xi$, was studied by Garc\'{\i}a-Mel\'{\i}an, Sabina de Lis and Sanabria-Garc\'{\i}a in \cite{GSS} and also by de la Calle, Sabina de Lis and Segura de Le\'{o}n in \cite{DSS}.

As far as we know, the best result concerning the existence of positive weak
solutions to the problem (\ref{Pw}), where the drift coefficient $a$ is
allowed to change sign, is the one obtained by Shouchuan and Papageorgiou in
\cite{SP}. In that paper, the authors combined the theory of nonlinear
operators of monotone type with comparison techniques to prove the existence
of $\lambda_{a}>0$ for which (\ref{Pw}) admits a weak solution $u_{a}\in
C_{0}^{1}(\overline{\Omega}),$ under the single hypothesis%
\begin{equation}
\left\Vert a\right\Vert _{\infty}<\lambda_{1,p}^{1/p}(\Omega).\label{HSP}%
\end{equation}
Here $\lambda_{1,p}(\Omega)$ denotes the first eigenvalue of the Dirichlet
$p$-Laplacian in $\Omega$ (see Section \ref{Basic}). 

Our main result in this paper improves upon theirs by eliminating hypothesis
(\ref{HSP}). It is stated as follows, where $R_{\Omega}$ denotes the
circumradius of $\Omega$ (i.e. the radius of the smallest ball containing
$\Omega$), $r_{\Omega}$ denotes the inradius of $\Omega$ (i.e. the radius of a
largest ball inscribed in $\Omega$), and the functions
\[
a^{+}:=\max\left\{  a,0\right\}  \text{ \ and \ }a^{-}:=\max\left\{
-a,0\right\}
\]
denote the positive and the negative parts of the function $a,$ respectively.

\begin{theorem}
\label{teomain}There exists a pair $(\lambda_{a},u_a)$ such that $\lambda_a >0$ and $u_{a}\in W_{0}^{1,p}(\Omega)\cap C^{1}(\overline
{\Omega})$ is a weak solution to the problem
\begin{equation}
\left\{
\begin{array}
[c]{lll}%
-\Delta_{p}u=\lambda_{a}\left\vert u\right\vert ^{p-2}u+a(x)\left\vert \nabla
u\right\vert ^{p-1} & \mathrm{in} & \Omega\\
u=0 & \mathrm{on} & \partial\Omega.
\end{array}
\right.  \label{Ppw}%
\end{equation}
Moreover, $u_a$ is positive in $\Omega$ and 
\begin{equation}
\left(  \frac{p}{p-1}\right)  ^{p-1}NR_{\Omega}^{-p}e^{-\left\Vert a^{+}\right\Vert
_{\infty}R_{\Omega}}\leq\lambda_{a}\leq e^{\left\Vert
a^{-}\right\Vert _{\infty}r_{\Omega}}\lambda_{1,p}\left(  B_{r_{\Omega}}%
(x_{0})\right)  ,\label{bounds}%
\end{equation}
where $B_{r_{\Omega}}(x_{0})\subset\Omega$ denotes a ball of radius
$r_{\Omega}$ centered at a point $x_{0}\in\Omega.$ If $a\geq0,$ then
$\lambda_{a}<\lambda_{1,p}(\Omega)$, and if $a\leq0,$ then $\lambda_{a}\geq\lambda
_{1,p}(\Omega).$
\end{theorem}

The particular form of the convective term in (\ref{Pw}) makes
the problem equivalent to a variational one when restricted to the positive,
radially symmetric weak solutions in a ball. Thus, working in this domain with the PDE
$-\Delta_{p}u=\left\vert u\right\vert ^{q-2}u+a(x)\left\vert \nabla
u\right\vert ^{p-1},$ for $1\leq q<p$, we first prove, by means of an
iterative method, an existence result when the coefficient $a$ has
radial symmetry. Then, we obtain a positive, radially
symmetric weak solution to the $1$-homogeneous problem by a
limiting procedure ($q\rightarrow p^{-}$).

We extend this approach for a general bounded domain $\Omega.$ In order to
obtain a positive weak solution $w_{q}\in W_{0}^{1,p}(\Omega)$ to the problem
\[
\left\{
\begin{array}
[c]{lll}%
-\Delta_{p}u=\left\vert u\right\vert ^{q-2}u+a(x)\left\vert \nabla
u\right\vert ^{p-1} & \mathrm{in} & \Omega\\
u=0 & \mathrm{on} & \partial\Omega,
\end{array}
\right.
\]
we apply a sub-supersolution result due to Cuesta Leon \cite{CL}. Our
subsolution is the zero extension to $\Omega$ of the positive, radially
symmetric weak solution to the problem
\[
\left\{
\begin{array}
[c]{lll}%
-\Delta_{p}u=\left\vert u\right\vert ^{q-2}u-\left\Vert a^{-}\right\Vert
_{\infty}\left\vert \nabla u\right\vert ^{p-1} & \mathrm{in} & B_{r_{\Omega}%
}(x_{0})\\
u=0 & \mathrm{on} & \partial B_{r_{\Omega}}(x_{0}).
\end{array}
\right.  \label{Pa-}%
\]
The supersolution is the restriction to $\Omega$ of the positive, radially
symmetric weak solution to the problem
\[
\left\{
\begin{array}
[c]{lll}%
-\Delta_{p}u=\left\vert u\right\vert ^{q-2}u+\left\Vert a^{+}\right\Vert
_{\infty}\left\vert \nabla u\right\vert ^{p-1} & \mathrm{in} & B_{R_{\Omega}%
}(x_{1})\\
u=0 & \mathrm{on} & \partial B_{R_{\Omega}}(x_{1}),
\end{array}
\right.  \label{Pa+}%
\]
where $B_{R_{\Omega}}(x_{1})$ is the smallest ball containing $\Omega.$

The essential tool for establishing the ordering between the subsolution and
the supersolution is the comparison principle by Abdellaoui and Peral
\cite{AP} (see Theorem \ref{Ape}) applied to $\vartheta_{q},$ the only
positive weak solution to the purely sublinear problem%
\[
\left\{
\begin{array}
[c]{lll}%
-\Delta_{p}u=\left\vert u\right\vert ^{q-2}u & \mathrm{in} & \Omega\\
u=0 & \mathrm{on} & \partial\Omega.
\end{array}
\right.
\]
This comparison principle guarantees that the subsolution and the supersolution lie below and above $\vartheta_{q},$ respectively.

Combining estimates derived from the sub-supersolutions method with a
well-known regularity result by Lieberman \cite{Lieb}, we obtain the pair
$\left(  \lambda_{a},u_{a}\right)  $ by letting $q\rightarrow p^{-},$ respectively
in
\[
\gamma_{q}:=\frac{1}{\left\Vert w_{q}\right\Vert _{\infty}^{q-p}}\text{ \ and
\ }u_{q}:=\frac{w_{q}}{\left\Vert w_{q}\right\Vert _{\infty}}.
\]

Our approach still produces the explicit bounds (\ref{bounds}) to $\lambda
_{a}$ that allow us to study the asymptotic behavior of $\lambda_{a}$ and
$u_{a},$ as $p\rightarrow\infty,$ considering that $a$ may depend on $p$, that
is, $a=a_{p}$. Thus, under the hypotheses%
\[
\limsup_{p\rightarrow\infty}\left\Vert a_{p}^{+}\right\Vert _{\infty}%
<r_{\Omega}^{-1}\ \ \text{and \ }\lim_{p\rightarrow\infty}\frac{\left\Vert
a_{p}^{-}\right\Vert _{\infty}}{p}=0,
\]
we show that
\[
\lim_{p\rightarrow\infty}\lambda_{a_{p}}^{1/p}=\lim_{p\rightarrow\infty}%
\frac{\left\Vert \nabla u_{a_{p}}\right\Vert _{p}}{\left\Vert u_{a_{p}%
}\right\Vert _{p}}=r_{\Omega}^{-1}.
\]
The homogeneity of (\ref{Ppw}) allows us to assume that $\left\Vert u_{a_{p}%
}\right\Vert _{p}=1.$ Thus, under this assumption we prove that, up to a
subsequence, $u_{a_{p}}$ converges uniformly to a nonnegative function
$u_{\infty}\in C_{0}(\overline{\Omega})\cap W^{1,\infty}(\Omega)$ such that
$\left\Vert u_{\infty}\right\Vert _{\infty}=1,$ $\left\Vert \nabla u_{\infty
}\right\Vert _{\infty}=r_{\Omega}^{-1},$ and 
$
0\leq u_{\infty}\leq r_{\Omega}^{-1}d_{\Omega},
$
where $d_{\Omega}\in C_{0}(\overline{\Omega})\cap W^{1,\infty}(\Omega)$
denotes the distance function to the boundary $\partial\Omega.$ These
properties imply that any maximum point of $u_{\infty}$ is also a maximum
point of $d_{\Omega}.$

This paper is organized as follows. In Section \ref{Basic} we introduce the
notation and recall some definitions and results that are used throughout the text.

In Section \ref{auxprob}, we present some properties of a variational problem
similar to (\ref{dot}). Such properties are used in Section \ref{secball} to
study problem (\ref{Pw}) when the domain is the ball $B_{R}(x_{0})$ centered at
$x_{0}\in\mathbb{R}^{N}$ with radius $R>0,$ and the coefficient is radially
symmetric, that is, $a(x)=a(\left\vert x-x_{0}\right\vert )$ with $a\in
C([0,R])$. We first establish, in Theorem \ref{teovqa}, the existence of a
positive, radially symmetric weak solution $w_{q,a}$ to the problem
\[
\left\{
\begin{array}
[c]{lll}%
-\Delta_{p}u=\left\vert u\right\vert ^{q-2}u+a\left(  \left\vert
x-x_{0}\right\vert \right)  \left\vert \nabla u\right\vert ^{p-1} &
\mathrm{in} & B_{R}(x_{0})\\
u=0 & \mathrm{on} & \partial B_{R}(x_{0}),
\end{array}
\right.
\]
for each $1\leq q<p,$ as well as estimates for $w_{q,a}$ and $\left\Vert \nabla
w_{q,a}\right\Vert _{\infty}.$ 
%The weak solutions to problems (\ref{Pa-}) and
%(\ref{Pa+}) are obtained by applying Theorem \ref{teovqa} to the particular
%case where the coefficient $a$ is constant.% 
The estimates derived in Theorem
\ref{teovqa} allow us to prove a more refined version of Theorem \ref{teomain}
for the particular case where the domain is a ball. In fact, in Theorem
\ref{teowa} we obtain, via a limiting argument, a positive, radially symmetric
weak solution $w_{a}$ to the problem
\[
\left\{
\begin{array}
[c]{lll}%
-\Delta_{p}w=\sigma_{a}\left\vert w\right\vert ^{p-2}w+a\left(  \left\vert
x-x_{0}\right\vert \right)  \left\vert \nabla w\right\vert ^{p-1} &
\mathrm{in} & B_{R}(x_{0})\\
w=0 & \mathrm{on} & \partial B_{R}(x_{0}),
\end{array}
\right.
\]
where
\[
\sigma_{a}:=\inf\left\{  \frac{%
%TCIMACRO{\dint \nolimits_{B_{R}(x_{0})}}%
%BeginExpansion
{\displaystyle\int\nolimits_{B_{R}(x_{0})}}
%EndExpansion
e^{-A(\left\vert x-x_{0}\right\vert )}\left\vert \nabla w\right\vert
^{p}\mathrm{d}x}{%
%TCIMACRO{\dint \nolimits_{B_{R}(x_{0})}}%
%BeginExpansion
{\displaystyle\int\nolimits_{B_{R}(x_{0})}}
%EndExpansion
e^{-A(\left\vert x-x_{0}\right\vert )}\left\vert w\right\vert ^{p}\mathrm{d}%
x},0\not \equiv w\in W_{0}^{1,p}(B_{R}(x_{0}))\right\}  ,
\]
and $A(r)=%
%TCIMACRO{\dint \nolimits_{0}^{r}}%
%BeginExpansion
{\displaystyle\int\nolimits_{0}^{r}}
%EndExpansion
a\left(  t\right)  \mathrm{d}t,$ for $r\in\left[  0,R\right]  .$ We also show
that $w_{a}$ is unique if $a$ does not change sign in $[0,R].$ In Section
\ref{secfinal} we apply the results developed in Section \ref{secball} to
prove Theorem \ref{teomain}. Section \ref{secasy} is dedicated to the
asymptotic behavior, as $p\rightarrow\infty,$ of $\lambda_{a_{p}}$ and
$u_{a_{p}}.$

For the sake of completeness we prove in the Appendix \ref{null} that if $\lambda \leq 0,$ then the only nonnegative weak solution to (\ref{Pw}) is the trivial one. 

\section{Definitions and known facts\label{Basic}}

In this section we introduce the notation and recall some definitions and results.

We will denote by $\left\Vert u\right\Vert _{r}$ the standard norm of a
function $u$ in the Lebesgue space $L^{r}(\Omega),$ $1\leq r\leq\infty.$ The
$L^{p}$-norm of $\left\vert \nabla u\right\vert ,$ the modulus of gradient of
a function $u\in W_{0}^{1,p}(\Omega),$ will be denoted by $\left\Vert
\nabla u\right\Vert _{p}.$  We recall that the Sobolev space $W_{0}%
^{1,p}(\Omega),$ $1<p<\infty,$ is the Banach space characterized as the
closure of $C_{c}^{\infty}(\Omega)$ with respect to the norm $\left\Vert
\nabla\cdot\right\Vert _{p}$, where $C_{c}^{\infty}(\Omega)$ denotes the class
of infinitely differentiable functions with compact support in $\Omega.$ Thus,
$\left\Vert \nabla u\right\Vert _{p}$ is the standard norm of a function $u\in
W_{0}^{1,p}(\Omega).$ It is well known that the embedding $W_{0}^{1,p}%
(\Omega)\hookrightarrow L^{p}(\Omega)$ is continuous and compact. Thus, 
\[
\lambda_{1,p}(\Omega):=\inf\left\{  \frac{\left\Vert \nabla u\right\Vert
_{p}^{p}}{\left\Vert u\right\Vert _{p}^{p}}:0\not \equiv u\in W_{0}%
^{1,p}(\Omega)\right\}>0  .
\]
We recall that $\lambda_{1,p}(\Omega)$ is also the first eigenvalue
of the Dirichlet $p$-Laplacian, that is, the least value of $\lambda>0$ such
that the Dirichlet problem
\[
\left\{
\begin{array}
[c]{lll}%
-\Delta_{p}u=\lambda\left\vert u\right\vert ^{p-2}u & \mathrm{in} & \Omega\\
u=0 & \mathrm{on} & \partial\Omega,
\end{array}
\right.
\]
has a nontrivial weak solution.

Let us consider the Dirichlet problem
\begin{equation}
\left\{
\begin{array}
[c]{lll}%
-\Delta_{p}u=f(x,u,\nabla u) & \mathrm{in} & \Omega\\
u=0 & \mathrm{on} & \partial\Omega,
\end{array}
\right.  \label{dir}%
\end{equation}
where $f:\Omega\times\mathbb{R}\times\mathbb{R}^{N}$ is a Carath\'{e}odory
function and $\Omega$ is a bounded domain of $\mathbb{R}^{N}.$

A function $u\in W_{0}^{1,p}(\Omega)$ is a weak solution to (\ref{dir})
if
\[
{\int_{\Omega}|\nabla u|^{p-2}\nabla u\cdot\nabla\phi\,\mathrm{d}%
x=\int_{\Omega}f(x,u,\nabla u)\phi\,\mathrm{d}x,}\text{ for all }\phi\in
W_{0}^{1,p}(\Omega).
\]

The following global regularity result is the version for the Dirichlet
problem (\ref{dir}) of a well known result by Lieberman \cite[Theorem 1]{Lieb}.

\begin{theorem}
\label{Lie}Let $\Omega$ be a bounded domain of $\mathbb{R}^{N}$ with
$C^{1,\alpha}$ boundary. Suppose that $u\in W_{0}^{1,p}(\Omega)$ is a bounded
weak solution to the Dirichlet problem (\ref{dir}), where the Carath\'{e}odory
function $f$ is such that%
\begin{equation}
\left\vert f(x,s,\xi)\right\vert \leq K\left(  1+\left\vert \xi\right\vert
^{p}\right),  \label{hipf}%
\end{equation}
for all $(x,s,\xi)\in\Omega\times\lbrack-M,M]\times\mathbb{R}^{N},$ for
positive constants $K$ and $M.$ If $\left\Vert u\right\Vert _{\infty}\leq M,$
then there exists $0<\tau<1,$ depending only on $\alpha,$ $K,$ $p$ and $N,$
such that $u\in C^{1,\tau}(\overline{\Omega});$ moreover, $\left\Vert
u\right\Vert _{C^{1,\tau}}\leq C$ where $C$ is a positive constant that
depends only on $\alpha,$ $K,$ $p,$ $M,$ $N$ and $\Omega.$
\end{theorem}

By a subsolution to (\ref{dir}) we mean a function $u\in W^{1,p}(\Omega)\cap
C(\overline{\Omega})$ such that $u\leq0$ on $\partial\Omega$ and
\[
-\Delta_{p}u\leq f(x,u,\nabla u)\text{ in }\Omega.
\]
This inequality means that
\[
{\int_{\Omega}|\nabla u|^{p-2}\nabla u\cdot\nabla\phi\,\mathrm{d}x\leq
\int_{\Omega}f(x,u,\nabla u)\phi\,\mathrm{d}x,}%
\]
for all ${\phi}\in C_{c}^{\infty}(\Omega),$ ${\phi\geq0}.$

A supersolution to (\ref{dir}) is a function $u\in W^{1,p}(\Omega)\cap
C(\overline{\Omega})$ such that $u\geq0$ on $\partial\Omega$ and
\[
-\Delta_{p}u\geq f(x,u,\nabla u)\text{ in }\Omega,
\]
meaning that
\[
{\int_{\Omega}|\nabla u|^{p-2}\nabla u\cdot\nabla\phi\,\mathrm{d}x\geq
\int_{\Omega}f(x,u,\nabla u)\phi\,\mathrm{d}x,}%
\]
for all ${\phi}\in C_{c}^{\infty}(\Omega),$ ${\phi\geq0}.$

These notions of subsolution and supersolution are stronger than the standard
ones, but they are sufficient for our purposes in this paper.

We say that $\underline{u}$ and $\overline{u}$ form an ordered pair of
sub-supersolutions to (\ref{dir}) if $\underline{u}$ is a subsolution to
(\ref{dir}), $\overline{u}$ is a supersolution to (\ref{dir}), and
$\underline{u}\leq\overline{u}$ in $\Omega.$

Next, we present a simplified version of a result due to Cuesta Leon
\cite[Theorem 2.2]{CL}.

\begin{theorem}
\label{Leon}Let $\underline{u},\overline{u}\in W_{0}^{1,p}(\Omega)$ be an
ordered pair of sub-supersolutions to the Dirichlet problem (\ref{dir}).
Suppose that
\begin{equation}
\left\vert f(x,s,\xi)\right\vert \leq K_{1}+K_{2}\left\vert \xi\right\vert
^{r} \label{hipfgrad}%
\end{equation}
a.e. $x\in\Omega,$ for all $s\in\lbrack\underline{u}(x),\overline{u}(x)]$ and
$\xi\in\mathbb{R}^{N},$ where $K_{1}$ and $K_{2}$ are positive constants and
\[
0\leq r<p-1+\min\left\{  \frac{p}{N},1\right\}  .
\]
Then, (\ref{dir}) has at least one weak solution $u\in W_{0}^{1,p}(\Omega)$
such that $\underline{u}\leq u\leq\overline{u}$ a.e. in $\Omega.$
\end{theorem}

As part of our approach to prove Theorem %
\ref{teomain} in Section \ref{secfinal}, we will apply Theorem \ref{Leon} to
the function 
\[
f(x,s,\xi )=\left\vert s\right\vert ^{q-2}s+a(x)\left\vert \xi \right\vert
^{p-1},
\]%
where $1\leq q<p.$ Alternatively, we could apply a result due to Boccardo,
Murat and Puel \cite[Theorem 2.1]{BMP} established for a Carath\'{e}odory function 
$f$ satisfying 
\begin{equation}
\left\vert f(x,s,\xi )\right\vert \leq C(\left\vert s\right\vert )\left(
1+\left\vert \xi \right\vert ^{p}\right) ,  \label{Boc}
\end{equation}%
for an increasing function $C:\mathbb{R}^{+}\rightarrow \mathbb{R}^{+}.$ Our
choice is merely for convenience, since the result by Cuesta Leon is for the $%
p$-Laplacian while the result by Boccardo, Murat and Puel is for a more
general operator.

It is well known that for each $q\in\lbrack1,p)$ the Dirichlet problem
\begin{equation}
\left\{
\begin{array}
[c]{lll}%
-\Delta_{p}v=\left\vert v\right\vert ^{q-2}v & \mathrm{in} & \Omega\\
v=0 & \mathrm{on} & \partial\Omega
\end{array}
\right.  \label{vq}%
\end{equation}
has a unique positive weak solution $\vartheta_{q}\in W_{0}^{1,p}(\Omega)\cap
L^{\infty}(\Omega).$ Hence, according to Theorem \ref{Lie}, $\vartheta_{q}\in
C^{1,\tau}(\overline{\Omega})$ for some $\tau\in\left(  0,1\right)  .$
The next result can be found in \cite[Theorem 8]{BJEM}.

\begin{theorem}
\label{sublin}Let $\Omega\subset\mathbb{R}^{N}$ be a bounded domain with a
Lipschitz boundary. One has%
\[
\lim_{q\rightarrow p^{-}}\frac{1}{\left\Vert \vartheta_{q}\right\Vert
_{\infty}^{p-q}}=\lambda_{1,p}(\Omega).
\]

\end{theorem}

The following comparison principle by Abdellaoui and Peral can be found in
\cite[Lemma 4.1]{AP}.

\begin{theorem}
\label{Ape}Let $\Omega\subset\mathbb{R}^{N}$ be a bounded domain with a
$C^{2}$ boundary. Let $f$ be a nonnegative, continuous function such that
$f(u)/u^{p-1}$ is decreasing. Assume that $u,v\in W_{0}^{1,p}(\Omega)\cap
C^{1}(\Omega)$ are such that%
\[
\left\{
\begin{array}
[c]{l}%
-\Delta_{p}u\geq f(u),\text{ \ }u>0\text{ \ in }\Omega\\
-\Delta_{p}v\leq f(v),\text{ \ }v>0\text{ \ in }\Omega.
\end{array}
\right.
\]
Then, $v\leq u$ in $\Omega.$
\end{theorem}

\begin{corollary}
\label{abpe2}Suppose that $\underline{v}\in W^{1,p}(\Omega)\cap C(\overline
{\Omega})$ is a nontrivial function satisfying
\[
-\Delta_{p}\underline{v}\leq\underline{v}^{q-1},\text{ \ }\underline{v}%
\geq0\ \ \text{in }\Omega.
\]
Then, $\underline{v}\leq\vartheta_{q}$ \ in $\Omega.$
\end{corollary}

\begin{proof}
Let $w\in W_{0}^{1,p}(\Omega)$ be such that $-\Delta_{p}w=\underline{v}%
^{q-1}.$ Theorem \ref{Lie} and a standard comparison principle guarantee that $w\in
W_{0}^{1,p}(\Omega)\cap C^{1}(\Omega)$ and $0\leq\underline{v}\leq w,$ respectively. Since $\underline{v}$ is nonnegative and nontrivial, we have $-\Delta_{p}w \geq 0$ and $w>0$ in $\Omega$ (according to the strong maximum principle). Thus, by applying Theorem \ref{Ape} we obtain
$w\leq\vartheta_{q}$ in $\Omega.$ Consequently, $\underline{v}\leq
w\leq\vartheta_{q}$ in $\Omega.$
\end{proof}

\section{An auxiliary variational problem\label{auxprob}}

In this section we fix an arbitrary $x_{0}\in\mathbb{R}^{N}$ and consider the
Dirichlet problem
\begin{equation}
\left\{
\begin{array}
[c]{lll}%
-\Delta_{p}u=\lambda\left\vert u\right\vert ^{q-2}u-a(\left\vert
x-x_{0}\right\vert )\left(  \dfrac{x-x_{0}}{\left\vert x-x_{0}\right\vert
}\cdot\nabla u\right)  \left\vert \nabla u\right\vert ^{p-2} & \mathrm{in} &
\Omega\\
u=0 & \mathrm{on} & \partial\Omega,
\end{array}
\right.  \label{axp}%
\end{equation}
where $a\in C([0,\infty))$ and $1\leq q\leq p$, with $p>1$.

A weak solution to this problem is a function $u\in W_{0}^{1,p}(\Omega)$ such
that
\begin{equation}%
%TCIMACRO{\dint \nolimits_{\Omega}}%
%BeginExpansion
{\displaystyle\int\nolimits_{\Omega}}
%EndExpansion
\left\vert \nabla u\right\vert ^{p-2}\nabla u\cdot\nabla\phi\mathrm{d}x=%
%TCIMACRO{\dint \nolimits_{\Omega}}%
%BeginExpansion
{\displaystyle\int\nolimits_{\Omega}}
%EndExpansion
\left[  \lambda\left\vert u\right\vert ^{q-2}u-a(\left\vert x-x_{0}\right\vert
)\left(  \dfrac{x-x_{0}}{\left\vert x-x_{0}\right\vert }\cdot\nabla u\right)
\left\vert \nabla u\right\vert ^{p-2}\right]  \phi\mathrm{d}x, \label{wpq}%
\end{equation}
for all $\phi\in W_{0}^{1,p}(\Omega).$

Let us set
\[
A(r):=%
%TCIMACRO{\dint \nolimits_{0}^{r}}%
%BeginExpansion
{\displaystyle\int\nolimits_{0}^{r}}
%EndExpansion
a(t)\mathrm{d}t.
\]
Testing (\ref{wpq}) with the function $\phi=e^{-A(\left\vert x-x_{0}%
\right\vert )}{\varphi}$, for some arbitrary ${\varphi\in W}_{0}^{1,p}%
(\Omega),$ one can verify that a weak solution to (\ref{axp}) is also a weak
solution to the problem
\begin{equation}
\left\{
\begin{array}
[c]{lll}%
-\operatorname{div}(e^{-A(\left\vert x-x_{0}\right\vert )}\left\vert \nabla
w\right\vert ^{p-2}\nabla w)=\lambda e^{-A(\left\vert x-x_{0}\right\vert
)}\left\vert w\right\vert ^{q-2}w & \mathrm{in} & \Omega\\
w=0 & \mathrm{on} & \partial\Omega.
\end{array}
\right.  \label{eqp}%
\end{equation}

A weak solution to this problem is a function $w\in W_{0}^{1,p}(\Omega)$ such
that
\begin{equation}%
%TCIMACRO{\dint \nolimits_{\Omega}}%
%BeginExpansion
{\displaystyle\int\nolimits_{\Omega}}
%EndExpansion
e^{-A(\left\vert x-x_{0}\right\vert )}\left\vert \nabla w\right\vert
^{p-2}\nabla w\cdot\nabla\phi\mathrm{d}x=\lambda%
%TCIMACRO{\dint \nolimits_{\Omega}}%
%BeginExpansion
{\displaystyle\int\nolimits_{\Omega}}
%EndExpansion
e^{-A(\left\vert x-x_{0}\right\vert )}\left\vert w\right\vert ^{q-2}%
w\phi\mathrm{d}x, \label{wepq}%
\end{equation}
for all $\phi\in W_{0}^{1,p}(\Omega).$

Testing (\ref{wepq}) with the function $\phi=e^{A(\left\vert x-x_{0}%
\right\vert )}{\varphi}$, for some arbitrary ${\varphi\in W}_{0}^{1,p}%
(\Omega),$ one can check that a weak solution to (\ref{eqp}) is also a weak
solution to (\ref{axp}). Therefore, (\ref{axp}) and (\ref{eqp}) are equivalent problems.

Problem (\ref{eqp}) is variational: its weak solutions are critical points
of the energy functional 
\[
E_{q}\left(  w\right)  :=\frac{1}{p}%
%TCIMACRO{\dint \nolimits_{\Omega}}%
%BeginExpansion
{\displaystyle\int\nolimits_{\Omega}}
%EndExpansion
e^{-A(\left\vert x-x_{0}\right\vert )}\left\vert \nabla w\right\vert
^{p}\mathrm{d}x-\frac{\lambda}{q}%
%TCIMACRO{\dint \nolimits_{\Omega}}%
%BeginExpansion
{\displaystyle\int\nolimits_{\Omega}}
%EndExpansion
e^{-A(\left\vert x-x_{0}\right\vert )}\left\vert w\right\vert ^{q}%
\mathrm{d}x, \textit{  }u\in W_{0}^{1,p}(\Omega).
\]

Taking into account that the function $e^{A(\left\vert x-x_{0}\right\vert )}$ is bounded above and below by positive constants, 
%\[
%0<\min_{x\in\overline{\Omega}}e^{A(\left\vert x-x_{0}\right\vert )}\leq
%\max_{x\in\overline{\Omega}}e^{A(\left\vert x-x_{0}\right\vert )}<\infty
%\]
one can verify that
\[
W_{0}^{1,p}(\Omega)\ni w\mapsto%
%TCIMACRO{\dint \nolimits_{\Omega}}%
%BeginExpansion
\left({\displaystyle \int\nolimits_{\Omega}}
%EndExpansion
e^{-A(\left\vert x-x_{0}\right\vert )}\left\vert \nabla w\right\vert
^{p}\mathrm{d}x\right)^{1/p}
\]
is a norm on the Sobolev space $W_{0}^{1,p}(\Omega)$ that is equivalent to the
standard norm and preserves the uniform convexity of the space.

The following result can be proved as in the well-known case $a\equiv0.$

\begin{proposition}
\label{lemwqa}Let $q\in\lbrack1,p)$ and $\lambda>0$ be fixed. The problem
\[
\left\{
\begin{array}
[c]{lll}%
-\operatorname{div}(e^{-A(\left\vert x-x_{0}\right\vert )}\left\vert \nabla
w\right\vert ^{p-2}\nabla w)=\lambda e^{-A(\left\vert x-x_{0}\right\vert
)}\left\vert w\right\vert ^{q-2}w & \mathrm{in} & \Omega\\
w=0 & \mathrm{on} & \partial\Omega
\end{array}
\right.
\]
has a unique weak solution $w_{q,a}\in W_{0}^{1,p}(\Omega)\cap C^{1}%
(\overline{\Omega})$ that is positive in $\Omega.$ Moreover,
$E_{q}\left(  w_{q,a}\right)  \leq E_{q}\left(  w\right)  $ for all $w\in
W_{0}^{1,p}(\Omega).$
\end{proposition}

Now, let us define
\[
\sigma_{a}(\Omega):=\inf\left\{  Q_{a}(w):w\in W_{0}^{1,p}(\Omega
)\setminus\left\{  0\right\}  \right\}  ,
\]
where%
\[
Q_{a}(w):=\left(
%TCIMACRO{\dint \nolimits_{\Omega}}%
%BeginExpansion
{\displaystyle\int\nolimits_{\Omega}}
%EndExpansion
e^{-A(\left\vert x-x_{0}\right\vert )}\left\vert \nabla w\right\vert
^{p}\mathrm{d}x\right)  \left(
%TCIMACRO{\dint \nolimits_{\Omega}}%
%BeginExpansion
{\displaystyle\int\nolimits_{\Omega}}
%EndExpansion
e^{-A(\left\vert x-x_{0}\right\vert )}\left\vert w\right\vert ^{p}%
\mathrm{d}x\right)  ^{-1}.
\]

The following proposition can also be proved as in the well-known case
$a\equiv0.$

\begin{proposition}
\label{lemwpa}There exists $w_{a}\in W_{0}^{1,p}(\Omega)\cap C^{1}%
(\overline{\Omega})$ such that $w_{a}>0$ in $\Omega,$
\[
\left\Vert w_{a}\right\Vert _{\infty
}=1\text{ and   }\sigma_{a}(\Omega)=Q_{a}(w_{a}).
\]
Moreover:

\begin{enumerate}
\item if $Q_{a}(w)=\sigma_{a}(\Omega),$ then $w=kw_{a}$ for some constant
$k\not =0;$

\item the function $w_{a}$ is a positive weak solution to the eigenvalue
problem
\begin{equation}
\left\{
\begin{array}
[c]{lll}%
-\operatorname{div}(e^{-A(\left\vert x-x_{0}\right\vert )}\left\vert \nabla
w\right\vert ^{p-2}\nabla w)=\lambda e^{-A(\left\vert x-x_{0}\right\vert
)}\left\vert w\right\vert ^{p-2}w & \mathrm{in} & \Omega\\
w=0 & \mathrm{on} & \partial\Omega,
\end{array}
\right.  \label{aux}%
\end{equation}
with $\lambda=\sigma_{a}(\Omega);$

\item if a function $\widetilde{w}\in W_{0}^{1,p}(\Omega)\setminus\left\{
0\right\}  $ is a positive weak solution to (\ref{aux}), then $\lambda
=\sigma_{a}(\Omega)$ and $\widetilde{w}=cw_{a}$ for some constant $c>0.$
\end{enumerate}
\end{proposition}

\begin{remark}
\label{center}It is clear that if $a\equiv0,$ then $\sigma_{0}(\Omega
)=\lambda_{1,p}(\Omega)$ and $w_{0}$ is the positive and $L^{\infty}%
$-normalized eigenfunction corresponding $\lambda_{1,p}(\Omega).$
\end{remark}

\begin{remark}
Problem (\ref{axp}), in the case where $q=p,$ takes the form of the eigenvalue
problem (\ref{B}) for the function
\[
B(x,\xi)=a(\left\vert x-x_{0}\right\vert )\left(  \frac{x-x_{0}}{\left\vert
x-x_{0}\right\vert }\cdot\xi\right)  \left\vert \xi\right\vert ^{p-2}.
\]
If either $x_{0}\notin \overline{\Omega}$, or $x_{0}\in \overline{\Omega}$ and $a(0)=0$, then $B$ is continuous on $\overline{\Omega}\times\mathbb{R}^{N}$ and, therefore, Proposition \ref{lemwpa} follows from \cite[Theorem 1.1]{GST}. However, the equivalence between problems \eqref{axp}, with $q=p$, and \eqref{aux} requires no additional assumptions on $x_{0}$ or $a(0)$. In fact, the boundedness of $a$ is sufficient to address problem \eqref{aux} by standard variational methods and establish Proposition \ref{lemwpa}. This is relevant to our approach, since in the subsequent sections we will deal with balls centered at $x_{0}$.

\end{remark}

In the next proposition, we obtain estimates for $\lambda>0$ for which the
problem
\begin{equation}
\left\{
\begin{array}
[c]{lll}%
-\Delta_{p}w=\lambda\left\vert w\right\vert ^{p-2}w+a(\left\vert
x-x_{0}\right\vert )\left\vert \nabla w\right\vert ^{p-1} & \mathrm{in} &
\Omega\\
w=0 & \mathrm{on} & \partial\Omega,
\end{array}
\right.  \label{aux0}%
\end{equation}
admits a positive weak solution in $W_{0}^{1.p}(\Omega)\cap C^{1}%
(\overline{\Omega}).$ We are only able to analyse the case in
which $a$ has a definite sign. In this case, our proof shows that if
$\lambda=\sigma_{a}(\Omega),$ then (\ref{aux0}) can admit at most one positive
weak solution. Furthermore, if such a solution exists, then it must be a
multiple of the minimizer $w_{a}$ corresponding to $\sigma_{a}(\Omega).$ Since
the EDPs satisfied by $w$ and $w_{a}$ are different from one another, this
conclusion suggests that
\[
\left\vert \nabla w\right\vert =\pm\frac{x-x_{0}}{\left\vert x-x_{0}%
\right\vert }\cdot\nabla w
\]
and points to the conjecture that $\Omega$ must be a ball centered at $x_{0}$.
This is an interesting question to investigate, as it may lead to a nonexistence
result when $\Omega$ is not such a ball.

\begin{proposition}
\label{nsym}Suppose that $w\in W_{0}^{1.p}(\Omega)\cap C^{1}(\overline{\Omega
})$ is a positive weak solution to the problem (\ref{aux0}) for some
$\lambda>0.$ If $a\geq0,$ then $\lambda\leq\sigma_{a}(\Omega),$ and if
$a\leq0,$ then $\sigma_{a}(\Omega)\leq\lambda.$ In both cases, if
$\lambda=\sigma_{a}(\Omega),$ then $w=cw_{a}$ for some positive constant $c.$
\end{proposition}

\begin{proof}
Let us first assume that $a\geq0$ in $\Omega.$ Let $w_{a}\in W_{0}%
^{1,p}(\Omega)$ be the positive minimizer corresponding to $\sigma_{a}\left(
\Omega\right)  .$

Testing the weak formulation of the problem (\ref{aux0}) with
$e^{-A(\left\vert x-x_{0}\right\vert )}\phi,$ where $\phi\in W_{0}%
^{1,p}(\Omega),$ $\phi\geq0,$ we find%
\[%
%TCIMACRO{\dint \nolimits_{\Omega}}%
%BeginExpansion
{\displaystyle\int\nolimits_{\Omega}}
%EndExpansion
e^{-A(\left\vert x-x_{0}\right\vert )}\left\vert \nabla w\right\vert
^{p-2}\nabla w\cdot\nabla\phi\mathrm{d}x=\lambda%
%TCIMACRO{\dint \nolimits_{\Omega}}%
%BeginExpansion
{\displaystyle\int\nolimits_{\Omega}}
%EndExpansion
e^{-A(\left\vert x-x_{0}\right\vert )}w^{p-1}\phi\mathrm{d}x+I(\phi)
\]
where%
\[
I(\phi):=%
%TCIMACRO{\dint \nolimits_{\Omega}}%
%BeginExpansion
{\displaystyle\int\nolimits_{\Omega}}
%EndExpansion
e^{-A(\left\vert x-x_{0}\right\vert )}a(x-x_{0})\left(  \left\vert \nabla
w\right\vert +\frac{x-x_{0}}{\left\vert x-x_{0}\right\vert }\cdot\nabla
w\right)  \left\vert \nabla w\right\vert ^{p-2}\phi\mathrm{d}x.
\]
As $a\geq0$ and
\[
\left\vert \nabla w\right\vert +\frac{x-x_{0}}{\left\vert x-x_{0}\right\vert
}\cdot\nabla w\geq0,
\]
we have that $I(\phi)\geq0$. Consequently,%
\[%
%TCIMACRO{\dint \nolimits_{\Omega}}%
%BeginExpansion
{\displaystyle\int\nolimits_{\Omega}}
%EndExpansion
e^{-A(\left\vert x-x_{0}\right\vert )}\left\vert \nabla w\right\vert
^{p-2}\nabla w\cdot\nabla\phi\mathrm{d}x\geq\lambda%
%TCIMACRO{\dint \nolimits_{\Omega}}%
%BeginExpansion
{\displaystyle\int\nolimits_{\Omega}}
%EndExpansion
e^{-A(\left\vert x-x_{0}\right\vert )}w^{p-1}\phi\mathrm{d}x.
\]

Let $\epsilon>0$ and choose $\phi:=\frac{w_{a}^{p}}{(w+\epsilon)^{p-1}}.$ It
follows that%
\[%
%TCIMACRO{\dint \nolimits_{\Omega}}%
%BeginExpansion
{\displaystyle\int\nolimits_{\Omega}}
%EndExpansion
e^{-A(\left\vert x-x_{0}\right\vert )}\left\vert \nabla w\right\vert
^{p-2}\nabla w\cdot\nabla\left(  \frac{w_{a}^{p}}{(w+\epsilon)^{p-1}}\right)
\mathrm{d}x\geq\lambda%
%TCIMACRO{\dint \nolimits_{\Omega}}%
%BeginExpansion
{\displaystyle\int\nolimits_{\Omega}}
%EndExpansion
e^{-A(\left\vert x-x_{0}\right\vert )}\left(  \frac{w}{w+\epsilon}\right)
^{p-1}w_{a}^{p}\mathrm{d}x.
\]

Now, recall that%
\[
\sigma_{a}(\Omega)%
%TCIMACRO{\dint \nolimits_{\Omega}}%
%BeginExpansion
{\displaystyle\int\nolimits_{\Omega}}
%EndExpansion
e^{-A(\left\vert x-x_{0}\right\vert )}w_{a}^{p}\mathrm{d}x=%
%TCIMACRO{\dint \nolimits_{\Omega}}%
%BeginExpansion
{\displaystyle\int\nolimits_{\Omega}}
%EndExpansion
e^{-A(\left\vert x-x_{0}\right\vert )}\left\vert \nabla w_{a}\right\vert
^{p}\mathrm{d}x.
\]
It follows from Picone's inequality (see \cite{AH}) that%
\begin{align*}
\left\vert \nabla w_{a}\right\vert ^{p}  &  \geq\left\vert \nabla
(w+\epsilon)\right\vert ^{p-2}\nabla(w+\epsilon)\cdot\nabla\left(  \frac
{w_{a}^{p}}{(w+\epsilon)^{p-1}}\right) \\
&  =\left\vert \nabla w\right\vert ^{p-2}\nabla w\cdot\nabla\left(
\frac{w_{a}^{p}}{(w+\epsilon)^{p-1}}\right)  .
\end{align*}
Hence,%
\[
\sigma_{a}(\Omega)%
%TCIMACRO{\dint \nolimits_{\Omega}}%
%BeginExpansion
{\displaystyle\int\nolimits_{\Omega}}
%EndExpansion
e^{-A(\left\vert x-x_{0}\right\vert )}w_{a}^{p}\mathrm{d}x\geq\lambda%
%TCIMACRO{\dint \nolimits_{\Omega}}%
%BeginExpansion
{\displaystyle\int\nolimits_{\Omega}}
%EndExpansion
e^{-A(\left\vert x-x_{0}\right\vert )}\left(  \frac{w}{w+\epsilon}\right)
^{p-1}w_{a}^{p}\mathrm{d}x.
\]
Letting $\epsilon\rightarrow0^{+}$ we find $\sigma_{a}(\Omega)\geq\lambda.$ If
$\lambda=\sigma_{a}(\Omega),$ then we have the inequality%
\[%
%TCIMACRO{\dint \nolimits_{\Omega}}%
%BeginExpansion
{\displaystyle\int\nolimits_{\Omega}}
%EndExpansion
e^{-A(\left\vert x-x_{0}\right\vert )}\left[  \left\vert \nabla w_{a}%
\right\vert ^{p}-\left\vert \nabla w\right\vert ^{p-2}\nabla w\cdot
\nabla\left(  \frac{w_{a}^{p}}{w^{p-1}}\right)  \right]  \mathrm{d}x\leq 0.
\]
As the integrand is nonnegative, we conclude that it is zero $a.e.$ in
$\Omega$ and this implies that $w=cw_{a}$ for some positive constant $c,$
according to \cite{AH}.

If $a\leq0,$ the proof is similar. Testing the weak formulation of the problem
(\ref{aux0}) with $e^{-A(\left\vert x-x_{0}\right\vert )}w$ we arrive at the
inequality
\[%
%TCIMACRO{\dint \nolimits_{\Omega}}%
%BeginExpansion
{\displaystyle\int\nolimits_{\Omega}}
%EndExpansion
e^{-A(\left\vert x-x_{0}\right\vert )}\left\vert \nabla w\right\vert
^{p}\mathrm{d}x\leq\lambda%
%TCIMACRO{\dint \nolimits_{\Omega}}%
%BeginExpansion
{\displaystyle\int\nolimits_{\Omega}}
%EndExpansion
e^{-A(\left\vert x-x_{0}\right\vert )}w^{p}\mathrm{d}x.
\]
Hence, by Picone's inequality we obtain%
\begin{align*}%
%TCIMACRO{\dint \nolimits_{\Omega}}%
%BeginExpansion
{\displaystyle\int\nolimits_{\Omega}}
%EndExpansion
e^{-A(\left\vert x-x_{0}\right\vert )}\left\vert \nabla w\right\vert
^{p}\mathrm{d}x  &  \geq%
%TCIMACRO{\dint \nolimits_{\Omega}}%
%BeginExpansion
{\displaystyle\int\nolimits_{\Omega}}
%EndExpansion
e^{-A(\left\vert x-x_{0}\right\vert )}\left\vert \nabla w_{a}\right\vert
^{p-2}\nabla w_{a}\cdot\nabla\left(  \frac{w^{p}}{(w_{a}+\epsilon)^{p-1}%
}\right)  \mathrm{d}x\\
&  =\sigma_{a}(\Omega)%
%TCIMACRO{\dint \nolimits_{\Omega}}%
%BeginExpansion
{\displaystyle\int\nolimits_{\Omega}}
%EndExpansion
e^{-A(\left\vert x-x_{0}\right\vert )}\left(  \frac{w_{a}}{w_{a}+\epsilon
}\right)  ^{p-1}w^{p}\mathrm{d}x,
\end{align*}
so that (after letting $\epsilon\rightarrow0^{+}$),%
\[
\lambda%
%TCIMACRO{\dint \nolimits_{\Omega}}%
%BeginExpansion
{\displaystyle\int\nolimits_{\Omega}}
%EndExpansion
e^{-A(\left\vert x-x_{0}\right\vert )}w^{p}\mathrm{d}x\geq\sigma_{a}(\Omega)%
%TCIMACRO{\dint \nolimits_{\Omega}}%
%BeginExpansion
{\displaystyle\int\nolimits_{\Omega}}
%EndExpansion
e^{-A(\left\vert x-x_{0}\right\vert )}w^{p}\mathrm{d}x.
\]
It follows that $\lambda\geq\sigma_{a}(\Omega),$ with the equality occurring
only if $w=cw_{a}$ for some positive constant $c.$
\end{proof}

\section{Radial solutions in a ball\label{secball}}

In this section,%
\[
B_{R}(x_{0}):=\left\{  x\in\mathbb{R}^{N}:\left\vert x-x_{0}\right\vert
<R\right\}  ,
\]
denotes the ball centered at $x_{0}\in\mathbb{R}^{N}$ with radius $R>0,$ and
$a:[0,R]\rightarrow\mathbb{R}$ is a continuous function with primitive%
\[
A(r):=%
%TCIMACRO{\dint \nolimits_{0}^{r}}%
%BeginExpansion
{\displaystyle\int\nolimits_{0}^{r}}
%EndExpansion
a(t)\mathrm{d}t,\text{ \ }0\leq r\leq R.
\]

Let us define the constants%
\[
\alpha:=\min_{0\leq s\leq t\leq R}e^{A(t)-A(s)}\text{ \ and \ }\beta
:=\max_{0\leq s\leq t\leq R}e^{A(t)-A(s)}.
\]
Of course, $0<\alpha<\beta<\infty.$

We establish two existence results. The first one is stated in
the sequel, where $1\leq q<p,$ $\vartheta_{q}$ is the positive weak solution
to (\ref{vq}) in $B_{R}(x_{0})$ and $w_{q,a}$ is the only positive weak
solution to the auxiliary problem%
\begin{equation}
\left\{
\begin{array}
[c]{lll}%
-\operatorname{div}\left(  e^{-A(\left\vert x-x_{0}\right\vert )}\left\vert
\nabla w\right\vert ^{p-2}\nabla w\right)  =e^{-A(\left\vert x-x_{0}%
\right\vert )}\left\vert w\right\vert ^{q-2}w & \mathrm{in} & B_{R}(x_{0})\\
w=0 & \mathrm{on} & \partial B_{R}(x_{0}).
\end{array}
\right.  \label{eRq}%
\end{equation}

\begin{theorem}
\label{teovqa}For each $q\in\lbrack1,p)$ the function $w_{q,a}$ is radially
symmetric, radially decreasing, belongs to $C^{2}\left(  B_{R}(x_{0}%
)\setminus\left\{  x_{0}\right\}  \right)  \cap C^{1}\left(  \overline
{B_{R}(x_{0})}\right)  ,$ and is also the only positive, radially symmetric
weak solution to the problem
\begin{equation}
\left\{
\begin{array}
[c]{lll}%
-\Delta_{p}w=\left\vert w\right\vert ^{q-2}w+a(\left\vert x-x_{0}\right\vert
)\left\vert \nabla w\right\vert ^{p-1} & \mathrm{in} & B_{R}(x_{0})\\
w=0 & \mathrm{on} & \partial B_{R}(x_{0}).
\end{array}
\right.  \label{Rq}%
\end{equation}
Moreover,
\begin{equation}
\alpha^{\frac{1}{p-q}}\vartheta_{q}\leq w_{q,a}\leq\left(  \frac{p-1}%
{p}\right)  ^{\frac{p-1}{p-q}}\left(  \frac{\beta R^{p}}{N}\right)  ^{\frac
{1}{p-q}}\text{ in }B_{R}(x_{0}), \label{estvqa1}%
\end{equation}
and%
\begin{equation}
\left\Vert \nabla w_{q,a}\right\Vert _{\infty}\leq\left(  \frac{p-1}%
{p}\right)  ^{\frac{p-1}{p-q}}\left(  \frac{\beta R^{q}}{N}\right)  ^{\frac
{1}{p-q}}. \label{estvqa2}%
\end{equation}

\end{theorem}

Next we state the second existence result, where $w_{a}$ is the unique
positive, $L^{\infty}$-normalized minimizer corresponding to $\sigma
_{a}\left(  B_{R}(x_{0})\right)  .$

\begin{theorem}
\label{teowa}The following estimates hold:
\begin{equation}
\left(  \frac{p}{p-1}\right)  ^{p-1}\frac{N}{\beta R^{p}}\leq\sigma_{a}\left(
B_{R}(x_{0})\right)  \leq\frac{\lambda_{1,p}\left(  B_{R}(x_{0})\right)
}{\alpha}\label{estsa}%
\end{equation}
and
\begin{equation}
\left\Vert \nabla w_{a}\right\Vert _{\infty}\leq\left(  \frac{\beta R}%
{N}\sigma_{a}\left(  B_{R}(x_{0})\right)  \right)  ^{\frac{1}{p-1}%
}.\label{estwa}%
\end{equation}
Moreover, $w_{a}$ is a radially symmetric weak solution to the
$1$-homogeneous problem
\begin{equation}
\left\{
\begin{array}
[c]{lll}%
-\Delta_{p}w=\sigma_{a}\left(  B_{R}(x_{0})\right)  \left\vert w\right\vert
^{p-2}w+a(\left\vert x-x_{0}\right\vert )\left\vert \nabla w\right\vert ^{p-1}
& \mathrm{in} & B_{R}(x_{0})\\
w=0 & \mathrm{on} & \partial B_{R}(x_{0}).
\end{array}
\right.  \label{PBRa}%
\end{equation}
In addition, if $\widetilde{w}$ is any other positive, radially symmetric weak
solution to the problem
\begin{equation}
\left\{
\begin{array}
[c]{lll}%
-\Delta_{p}w=\lambda\left\vert w\right\vert ^{p-2}w+a(\left\vert
x-x_{0}\right\vert )\left\vert \nabla w\right\vert ^{p-1} & \mathrm{in} &
B_{R}(x_{0})\\
w=0 & \mathrm{on} & \partial B_{R}(x_{0}),
\end{array}
\right.  \label{PBRl}%
\end{equation}
for some $\lambda>0,$ then $\widetilde{w}=cw_{a},$ for some positive constant
$c,$ and $\lambda=\sigma_{a}\left(  B_{R}(x_{0})\right)  .$
\end{theorem}

Theorems \ref{teovqa} and \ref{teowa} are proved in Subsection
\ref{teosvw}.

\begin{remark}
If $a$ does not change sign in $[0,R],$ then positive multiples of $w_{a}$ are
the only positive weak solutions to (\ref{PBRa}). This fact is a direct
consequence of Proposition \ref{nsym}. An interesting question to be
investigated is whether (\ref{PBRa}) admits a nonradial solution when $a$
changes sign.
\end{remark}

The following monotonicity
result is a consequence of the radial symmetry of $w_{a}$.

\begin{corollary}
Suppose that $a_{1}(r)\leq a_{2}(r)$ for all $r\in\lbrack0,R].$ Then,
\[
\sigma_{a_{2}}(B_{R}(x_{0}))\leq\sigma_{a_{1}}(B_{R}(x_{0})).
\]

\end{corollary}

\begin{proof}
For $i\in\left\{  1,2\right\}  $ let $w_{a_{i}}\in W_{0}^{1.p}(B_{R}%
(x_{0}))\cap C^{1}(\overline{B_{R}(x_{0})})$ be such that
\[
\sigma_{a_{i}}(B_{R}(x_{0}))%
%TCIMACRO{\dint \nolimits_{B_{R}(x_{0})}}%
%BeginExpansion
{\displaystyle\int\nolimits_{B_{R}(x_{0})}}
%EndExpansion
e^{-A_{i}(\left\vert x-x_{0}\right\vert )}\left\vert w_{a_{i}}\right\vert
^{p}\mathrm{d}x=%
%TCIMACRO{\dint \nolimits_{B_{R}(x_{0})}}%
%BeginExpansion
{\displaystyle\int\nolimits_{B_{R}(x_{0})}}
%EndExpansion
e^{-A_{i}(\left\vert x-x_{0}\right\vert )}\left\vert \nabla w_{a_{i}%
}\right\vert ^{p}\mathrm{d}x
\]
where $A_{i}(r)=%
%TCIMACRO{\dint \nolimits_{0}^{r}}%
%BeginExpansion
{\displaystyle\int\nolimits_{0}^{r}}
%EndExpansion
a_{i}(t)\mathrm{d}t.$ Using
\[
\phi:=e^{-(A_{2}(\left\vert x-x_{0}\right\vert )-A_{1}(\left\vert
x-x_{0}\right\vert ))}w_{a_{1}}%
\]
as a test function in the weak formulation (\ref{wepq}) with $q=p,$
$\lambda=\sigma_{a_{1}}(B_{R}(x_{0}))$ and $A=A_{1},$ we find
\begin{equation}%
%TCIMACRO{\dint \nolimits_{B_{R}(x_{0})}}%
%BeginExpansion
{\displaystyle\int\nolimits_{B_{R}(x_{0})}}
%EndExpansion
e^{-A_{2}(\left\vert x-x_{0}\right\vert )}\left\vert \nabla w_{a_{1}%
}\right\vert ^{p}\mathrm{d}x+J(w_{a_{1}})=\sigma_{a_{1}}(B_{R}(x_{0}))%
%TCIMACRO{\dint \nolimits_{B_{R}(x_{0})}}%
%BeginExpansion
{\displaystyle\int\nolimits_{B_{R}(x_{0})}}
%EndExpansion
e^{-A_{2}(\left\vert x-x_{0}\right\vert )}\left\vert w_{a_{1}}\right\vert
^{p}\mathrm{d}x, \label{J<0}%
\end{equation}
where
\[
J(w_{a_{1}}):=%
%TCIMACRO{\dint \nolimits_{B_{R}(x_{0})}}%
%BeginExpansion
{\displaystyle\int\nolimits_{B_{R}(x_{0})}}
%EndExpansion
e^{-A_{2}(\left\vert x-x_{0}\right\vert )}(a_{2}(\left\vert x-x_{0}\right\vert
)-a_{1}(\left\vert x-x_{0}\right\vert ))\left\vert \nabla w_{a_{1}}\right\vert
^{p-1}w_{a_{1}}\mathrm{d}x.
\]
Here we have used the following consequence of the symmetry of $w_{a_{1}}$:%
\[
-\frac{x-x_{0}}{\left\vert x-x_{0}\right\vert }\cdot\nabla w_{a_{1}%
}=\left\vert \nabla w_{a_{1}}\right\vert .
\]
As $J(w_{a_{1}})\geq0,$ the minimizing property of $w_{a_{2}}$ and (\ref{J<0})
yield%
\[
\sigma_{a_{2}}(B_{R}(x_{0}))\leq\frac{%
%TCIMACRO{\dint \nolimits_{B_{R}(x_{0})}}%
%BeginExpansion
{\displaystyle\int\nolimits_{B_{R}(x_{0})}}
%EndExpansion
e^{-A_{2}(\left\vert x-x_{0}\right\vert )}\left\vert \nabla w_{a_{1}%
}\right\vert ^{p}\mathrm{d}x}{%
%TCIMACRO{\dint \nolimits_{B_{R}(x_{0})}}%
%BeginExpansion
{\displaystyle\int\nolimits_{B_{R}(x_{0})}}
%EndExpansion
e^{-A_{2}(\left\vert x-x_{0}\right\vert )}\left\vert w_{a_{1}}\right\vert
^{p}\mathrm{d}x}\leq\sigma_{a_{1}}(B_{R}(x_{0})).
\]

\end{proof}

\subsection{The ODE boundary value problem}

Let us consider the Dirichlet problems
\begin{equation}
\left\{
\begin{array}
[c]{lll}%
-\Delta_{p}w=f(w)+a(\left\vert x-x_{0}\right\vert )\left\vert \nabla
w\right\vert ^{p-1} & \mathrm{in} & B_{R}(x_{0})\\
w=0 & \mathrm{on} & \partial B_{R}(x_{0})
\end{array}
\right.  \label{Rf}%
\end{equation}
and%
\begin{equation}
\left\{
\begin{array}
[c]{lll}%
-\operatorname{div}\left(  e^{-A(\left\vert x-x_{0}\right\vert )}\left\vert
\nabla w\right\vert ^{p-2}\nabla w\right)  =e^{-A(\left\vert x-x_{0}%
\right\vert )}f(w) & \mathrm{in} & B_{R}(x_{0})\\
w=0 & \mathrm{on} & \partial B_{R}(x_{0}),
\end{array}
\right.  \label{eRf}%
\end{equation}
where $f:\mathbb{R}\rightarrow\mathbb{R}$ is a continuous function such that
$f(0)=0$ and $f(t)>0$ for all $t>0.$

With respect to positive radial solutions, both problems (\ref{Rf}) and
(\ref{eRf}) are equivalent, since they have the same radial version, nominally
the ODE boundary value problem
\begin{equation}
\left\{
\begin{array}
[c]{lll}%
\left(  e^{-A(r)}r^{N-1}\left\vert v^{\prime}\right\vert ^{p-1}\right)
^{\prime}=e^{-A(r)}r^{N-1}f(v) & \mathrm{if} & r\in(0,R)\\
v^{\prime}(0)=0=v(R). &  &
\end{array}
\right.  \label{erq}%
\end{equation}

Indeed, if $w(x)=\omega\left(  \left\vert x-x_{0}\right\vert \right)  $ is a
positive, radially symmetric solution to (\ref{Rf}), then the function
$\omega(r)$ must be a solution to the ordinary boundary value problem%
\begin{equation}
\left\{
\begin{array}
[c]{l}%
-\left(  r^{N-1}\left\vert v^{\prime}\right\vert ^{p-2}v^{\prime}\right)
^{\prime}=r^{N-1}f(v)+a(r)r^{N-1}\left\vert v^{\prime}\right\vert
^{p-1};\text{ }r\in(0,R)\\
v^{\prime}(0)=0=v(R).
\end{array}
\right.  \label{radq}%
\end{equation}
It is simple to verify that the positivity of $\omega(r)$ forces
$\omega^{\prime}(r)<0$ for all $r\in\lbrack0,R).$ Consequently, $\omega
^{\prime}=-\left\vert \omega^{\prime}\right\vert ,$ and
\[
\left(  r^{N-1}\left\vert \omega^{\prime}(r)\right\vert ^{p-1}\right)
^{\prime}-a(r)r^{N-1}\left\vert \omega^{\prime}(r)\right\vert ^{p-1}%
=r^{N-1}f(\omega(r)).
\]
Using the integrating factor $e^{-A(r)}$ we can see that this ODE is
equivalent to that in (\ref{erq}). Inversely, if $v$ is a positive solution to
(\ref{erq}), then $v^{\prime}=-\left\vert v^{\prime}\right\vert $ (note that
the ODE cannot have any other critical point besides $r=0$). Hence, after
eliminating the integrating factor $e^{-A(r)}$ and writing
\[
\left(  r^{N-1}\left\vert v^{\prime}(r)\right\vert ^{p-1}\right)  ^{\prime
}=-\left(  r^{N-1}\left\vert v^{\prime}(r)\right\vert ^{p-2}v^{\prime
}(r)\right)  ^{\prime},
\]
we arrive at the ODE in (\ref{radq}).

Through a formal integration, we can deduce the following Green-type formula for a nonnegative solution to
(\ref{erq})
\[
v(r)=\int_{r}^{R}\left(
%TCIMACRO{\dint \nolimits_{0}^{t}}%
%BeginExpansion
{\displaystyle\int\nolimits_{0}^{t}}
%EndExpansion
e^{A(t)-A(s)}\left(  \frac{s}{t}\right)  ^{N-1}f(v)\mathrm{d}s\right)
^{\frac{1}{p-1}}\mathrm{d}t,\text{ }0\leq r\leq R.
\]
This expression suggests that a positive solution to (\ref{erq}) can be
obtained as a fixed point of the operator $T:C\left(  [0,R]\right)
\rightarrow C\left(  [0,R]\right)  $, defined as%

\begin{equation}
\left(  Tv\right)  (r)=\int_{r}^{R}\left(
%TCIMACRO{\dint \nolimits_{0}^{t}}%
%BeginExpansion
{\displaystyle\int\nolimits_{0}^{t}}
%EndExpansion
e^{A(t)-A(s)}\left(  \frac{s}{t}\right)  ^{N-1}f(\left\vert v\right\vert
)\mathrm{d}s\right)  ^{\frac{1}{p-1}}\mathrm{d}t,\text{ }0\leq r\leq R.
\label{T}%
\end{equation}

Note that $\left(  Tv\right)  (R)=0$ and $Tv$ is of class $C^{1}$ in the
interval $(0,R],$ with%
\[
(Tv)^{\prime}(r)=-\left(  \frac{1}{r^{N-1}}%
%TCIMACRO{\dint \nolimits_{0}^{r}}%
%BeginExpansion
{\displaystyle\int\nolimits_{0}^{r}}
%EndExpansion
e^{A(r)-A(s)}s^{N-1}f(\left\vert v\right\vert )\mathrm{d}s\right)  ^{\frac
{1}{p-1}},\text{ }0<r\leq R.
\]
At $r=0$ the function $Tv$ is also of class $C^{1}$ with $(Tv)^{\prime}(0)=0,$
since%
\begin{align*}
(Tv)^{\prime}(0)  &  =\lim_{r\rightarrow0^{+}}\frac{(Tv)(r)-(Tv)(0)}{r}\\
&  =\lim_{r\rightarrow0^{+}}(Tv)^{\prime}(r)=-\left(  \lim_{r\rightarrow0^{+}%
}\frac{r^{N-1}f(\left\vert v(r)\right\vert )}{(N-1)r^{N-2}}\right)  ^{\frac
{1}{p-1}}=0.
\end{align*}
Therefore, $Tv\in C^{1}\left(  [0,R]\right)  .$

It is clear that $Tv\in C^{2}\left(  (0,R]\right)  .$ Moreover, the identity
\[
e^{-A(r)}r^{N-1}\left\vert (Tv)^{\prime}(r)\right\vert ^{p-1}=%
%TCIMACRO{\dint \nolimits_{0}^{r}}%
%BeginExpansion
{\displaystyle\int\nolimits_{0}^{r}}
%EndExpansion
e^{-A(s)}s^{N-1}f(\left\vert v(s)\right\vert )\mathrm{d}s,\text{ }0<r\leq R
\]
implies that $w=Tv\in C^{2}\left(  (0,R]\right)  \cap C^{1}\left(
[0,R]\right)  $ is a nonnegative function satisfying%
\begin{equation}
\left\{
\begin{array}
[c]{l}%
\left(  e^{-A(r)}r^{N-1}\left\vert w^{\prime}(r)\right\vert ^{p-1}\right)
^{\prime}=e^{-A(r)}r^{N-1}f(\left\vert v(r)\right\vert ),\text{ }0<r\leq R\\
w^{\prime}(0)=0=w(R).
\end{array}
\right.  \label{ax3}%
\end{equation}
Therefore, if $v\in C\left(  [0,R]\right)  $ is positive in $[0,R)$ and such
that $v=Tv,$ then $v\in C^{2}\left(  (0,R]\right)  \cap C^{1}\left(
[0,R]\right)  $ and $v$ is a solution to (\ref{erq}).

\begin{lemma}
\label{qw}Let $v\in C^{2}\left(  (0,R]\right)  \cap C^{1}\left(  [0,R]\right)
$ be a nonnegative solution to (\ref{erq}). The radially symmetric function
$w$ defined by
\[
w(x):=v(\left\vert x-x_{0}\right\vert ),\text{ \ }x\in B_{R}(x_{0})
\]
is a weak solution to both problems (\ref{Rf}) and (\ref{eRf}).
\end{lemma}

\begin{proof}
First we show that $w$ is a weak solution to (\ref{eRf}). Of course, $w\in
C^{2}\left(  B_{R}(x_{0})\setminus\left\{  x_{0}\right\}  \right)  \cap
C^{1}\left(  B_{R}(x_{0})\right)  $ and
\begin{equation}
-\operatorname{div}\left(  e^{-A\left(  \left\vert x-x_{0}\right\vert \right)
}\left\vert \nabla w(x)\right\vert ^{p-2}\nabla w(x)\right)  =e^{-A\left(
\left\vert x-x_{0}\right\vert \right)  }f(w(x))\label{div1}%
\end{equation}
for all $x\in B_{R}(x_{0})\setminus\left\{  x_{0}\right\}  .$

Let $\phi\in C_{c}^{\infty}\left(  B_{R}(x_{0})\right)  $ be fixed and write
\[%
%TCIMACRO{\dint \nolimits_{B_{R}(x_{0})}}%
%BeginExpansion
{\displaystyle\int\nolimits_{B_{R}(x_{0})}}
%EndExpansion
e^{-A\left(  \left\vert x-x_{0}\right\vert \right)  }\left\vert \nabla
w(x)\right\vert ^{p-2}\nabla w(x)\cdot\nabla\phi(x)\mathrm{d}x=I_{\epsilon
}+J_{\epsilon}%
\]
where $0<\epsilon<R,$%
\[
I_{\epsilon}:=%
%TCIMACRO{\dint \nolimits_{B_{\epsilon}(x_{0})}}%
%BeginExpansion
{\displaystyle\int\nolimits_{B_{\epsilon}(x_{0})}}
%EndExpansion
e^{-A\left(  \left\vert x-x_{0}\right\vert \right)  }\left\vert \nabla
w(x)\right\vert ^{p-2}\nabla w(x)\cdot\nabla\phi(x)\mathrm{d}x
\]
and%
\[
J_{\epsilon}:=%
%TCIMACRO{\dint \nolimits_{B_{R}(x_{0})\setminus B_{\epsilon}(x_{0})}}%
%BeginExpansion
{\displaystyle\int\nolimits_{B_{R}(x_{0})\setminus B_{\epsilon}(x_{0})}}
%EndExpansion
e^{-A\left(  \left\vert x-x_{0}\right\vert \right)  }\left\vert \nabla
w(x)\right\vert ^{p-2}\nabla w(x)\cdot\nabla\phi(x)\mathrm{d}x.
\]
We have that $\lim_{\epsilon\rightarrow0^{+}}I_{\epsilon}=0$ since
\[
\vert I_{\epsilon} \vert \leq\left\Vert \nabla w\right\Vert _{\infty}^{p-1}\left\Vert
\nabla\phi\right\Vert _{\infty}%
%TCIMACRO{\dint \nolimits_{B_{\epsilon}(x_{0})}}%
%BeginExpansion
{\displaystyle\int\nolimits_{B_{\epsilon}(x_{0})}}
%EndExpansion
e^{-A\left(  \left\vert x-x_{0}\right\vert \right)  }\mathrm{d}x
\]

As%
\[
\nabla w\left(  x\right)  =v^{\prime}(\left\vert x-x_{0}\right\vert
)\frac{x-x_{0}}{\left\vert x-x_{0}\right\vert }\text{ \ for all }x\in
B_{R}(x_{0})\setminus\left\{  x_{0}\right\}
\]
and $v^{\prime}\leq0,$ we that have%
\[
\left\vert \nabla w\left(  x\right)  \right\vert =-v^{\prime}(\left\vert
x-x_{0}\right\vert )
\]
and%
\begin{align*}
\left\vert \nabla w\left(  x\right)  \right\vert ^{p-2}\nabla w\left(
x\right)  \cdot\frac{x-x_{0}}{\left\vert x-x_{0}\right\vert }  &  =-\left\vert
v^{\prime}\left(  \left\vert x-x_{0}\right\vert \right)  \right\vert
^{p-1}\frac{x-x_{0}}{\left\vert x-x_{0}\right\vert }\cdot\frac{x-x_{0}%
}{\left\vert x-x_{0}\right\vert }\\
&  =-\left\vert v^{\prime}\left(  \left\vert x-x_{0}\right\vert \right)
\right\vert ^{p-1},
\end{align*}
for all $x\in B_{R}(x_{0})\setminus\left\{  x_{0}\right\}  .$

Thus, if $0<r\leq R,$ then
\[
\int_{\partial B_{r}(x_{0})}e^{-A\left(  \left\vert x-x_{0}\right\vert
\right)  }\left\vert \nabla w(x)\right\vert ^{p-2}\nabla w(x)\cdot
\frac{x-x_{0}}{\left\vert x-x_{0}\right\vert }\mathrm{d}S_{x}=-e^{-A(r)}%
\left\vert v^{\prime}(r)\right\vert ^{p-1}\int_{\partial B_{r}(x_{0})}%
\phi\left(  x\right)  \mathrm{d}S_{x}.
\]
Hence, taking into account that $%
%TCIMACRO{\dint \nolimits_{\partial B_{R}(x_{0})}}%
%BeginExpansion
{\displaystyle\int\nolimits_{\partial B_{R}(x_{0})}}
%EndExpansion
\phi\left(  x\right)  \mathrm{d}S_{x}=0,$ Divergence Theorem and (\ref{div1})
yield
\begin{align*}
J_{\epsilon}  &  =%
%TCIMACRO{\dint \nolimits_{B_{R}(x_{0})\setminus B_{\epsilon}(x_{0})}}%
%BeginExpansion
{\displaystyle\int\nolimits_{B_{R}(x_{0})\setminus B_{\epsilon}(x_{0})}}
%EndExpansion
\operatorname{div}\left(  \phi(x)e^{-A\left(  \left\vert x-x_{0}\right\vert
\right)  }\left\vert \nabla w(x)\right\vert ^{p-2}\nabla w(x)\right)
\mathrm{d}x\\
&  -%
%TCIMACRO{\dint \nolimits_{B_{R}(x_{0})\setminus B_{\epsilon}(x_{0})}}%
%BeginExpansion
{\displaystyle\int\nolimits_{B_{R}(x_{0})\setminus B_{\epsilon}(x_{0})}}
%EndExpansion
\operatorname{div}\left(  e^{-A\left(  \left\vert x-x_{0}\right\vert \right)
}\left\vert \nabla w(x)\right\vert ^{p-2}\nabla w(x)\right)  \phi
(x)\mathrm{d}x\\
&  =-e^{-A\left(  R\right)  }\left\vert v^{\prime}\left(  R\right)
\right\vert ^{p-1}%
%TCIMACRO{\dint \nolimits_{\partial B_{R}(x_{0})}}%
%BeginExpansion
{\displaystyle\int\nolimits_{\partial B_{R}(x_{0})}}
%EndExpansion
\phi\left(  x\right)  \mathrm{d}S_{x}+e^{-A\left(  \epsilon\right)
}\left\vert v^{\prime}\left(  \epsilon\right)  \right\vert ^{p-1}%
%TCIMACRO{\dint \nolimits_{\partial B_{\epsilon}(x_{0})}}%
%BeginExpansion
{\displaystyle\int\nolimits_{\partial B_{\epsilon}(x_{0})}}
%EndExpansion
\phi\left(  x\right)  \mathrm{d}S_{x}\\
&  +%
%TCIMACRO{\dint \nolimits_{B_{R}(x_{0})\setminus B_{\epsilon}(x_{0})}}%
%BeginExpansion
{\displaystyle\int\nolimits_{B_{R}(x_{0})\setminus B_{\epsilon}(x_{0})}}
%EndExpansion
e^{-A\left(  \left\vert x-x_{0}\right\vert \right)  }f(w(x))\phi
(x)\mathrm{d}x\\
&  =e^{-A\left(  \epsilon\right)  }\left\vert v^{\prime}\left(  \epsilon
\right)  \right\vert ^{p-1}%
%TCIMACRO{\dint \nolimits_{\partial B_{\epsilon}(x_{0})}}%
%BeginExpansion
{\displaystyle\int\nolimits_{\partial B_{\epsilon}(x_{0})}}
%EndExpansion
\phi\left(  x\right)  \mathrm{d}S_{x}+%
%TCIMACRO{\dint \nolimits_{B_{R}(x_{0})\setminus B_{\epsilon}(x_{0})}}%
%BeginExpansion
{\displaystyle\int\nolimits_{B_{R}(x_{0})\setminus B_{\epsilon}(x_{0})}}
%EndExpansion
e^{-A\left(  \left\vert x-x_{0}\right\vert \right)  }f(w(x))\phi
(x)\mathrm{d}x.
\end{align*}
Since
\[
0\leq e^{-A\left(  \epsilon\right)  }\left\vert v^{\prime}\left(
\epsilon\right)  \right\vert ^{p-1}\leq e^{-A\left(  \epsilon\right)
}\left\Vert \nabla w\right\Vert _{\infty}^{p-1}%
\]
we obtain%
\[
\lim_{\epsilon\rightarrow0^{+}}e^{-A\left(  \epsilon\right)  }\left\vert
v^{\prime}\left(  \epsilon\right)  \right\vert ^{p-1}%
%TCIMACRO{\dint \nolimits_{\partial B_{\epsilon}(x_{0})}}%
%BeginExpansion
{\displaystyle\int\nolimits_{\partial B_{\epsilon}(x_{0})}}
%EndExpansion
\phi\left(  x\right)  \mathrm{d}S_{x}=0,
\]
so that 
\[
\lim_{\epsilon\rightarrow0^{+}}J_{\epsilon}=%
%TCIMACRO{\dint \nolimits_{B_{R}(x_{0})}}%
%BeginExpansion
{\displaystyle\int\nolimits_{B_{R}(x_{0})}}
%EndExpansion
e^{-A\left(  \left\vert x-x_{0}\right\vert \right)  }f(w(x))\phi
(x)\mathrm{d}x.
\]
Consequently,
\begin{align*}%
%TCIMACRO{\dint \nolimits_{B_{R}(x_{0})}}%
%BeginExpansion
{\displaystyle\int\nolimits_{B_{R}(x_{0})}}
%EndExpansion
e^{-A\left(  \left\vert x-x_{0}\right\vert \right)  }\left\vert \nabla
w(x)\right\vert ^{p-2}\nabla w(x)\cdot\nabla\phi(x)\mathrm{d}x  &
=\lim_{\epsilon\rightarrow0^{+}}I_{\epsilon}+\lim_{\epsilon\rightarrow0^{+}%
}J_{\epsilon}\\
&  =%
%TCIMACRO{\dint \nolimits_{B_{R}(x_{0})}}%
%BeginExpansion
{\displaystyle\int\nolimits_{B_{R}(x_{0})}}
%EndExpansion
e^{-A\left(  \left\vert x-x_{0}\right\vert \right)  }f(w(x))\phi
(x)\mathrm{d}x.
\end{align*}
As $C_{c}^{\infty}\left(  B_{R}(x_{0})\right)  $ is dense in $W_{0}%
^{1,p}(B_{R}(x_{0}))$, we conclude that the identity%
\begin{equation}%
%TCIMACRO{\dint \nolimits_{B_{R}(x_{0})}}%
%BeginExpansion
{\displaystyle\int\nolimits_{B_{R}(x_{0})}}
%EndExpansion
e^{-A\left(  \left\vert x-x_{0}\right\vert \right)  }\left\vert \nabla
w\right\vert ^{p-2}\nabla w\cdot\nabla\phi\mathrm{d}x=%
%TCIMACRO{\dint \nolimits_{B_{R}(x_{0})}}%
%BeginExpansion
{\displaystyle\int\nolimits_{B_{R}(x_{0})}}
%EndExpansion
e^{-A\left(  \left\vert x-x_{0}\right\vert \right)  }f(w)\phi\mathrm{d}x
\label{wew}%
\end{equation}
holds for all $\phi\in W_{0}^{1,p}(B_{R}(x_{0})),$ that is, $w$ is a weak
solution to (\ref{eRf}).

Now, for an arbitrary ${\varphi}\in W_{0}^{1,p}(B_{R}(x_{0}))$ we use
$\phi=e^{A\left(  \left\vert x-x_{0}\right\vert \right)  }{\varphi}$ in
(\ref{wew}) and, after some simple calculations, taking into account that%
\begin{align*}
\left\vert \nabla w\right\vert ^{p-2}\nabla w\cdot\nabla e^{A\left(
\left\vert x-x_{0}\right\vert \right)  }  &  =\left\vert \nabla w\right\vert
^{p-2}\nabla w\cdot a(\left\vert x-x_{0}\right\vert )\frac{x-x_{0}}{\left\vert
x-x_{0}\right\vert }\\
&  =-a(\left\vert x-x_{0}\right\vert )\left\vert \nabla w\right\vert ^{p-1},
\end{align*}
we arrive at the equality
\[%
%TCIMACRO{\dint \nolimits_{B_{R}(x_{0})}}%
%BeginExpansion
{\displaystyle\int\nolimits_{B_{R}(x_{0})}}
%EndExpansion
\left\vert \nabla w\right\vert ^{p-2}\nabla w\cdot\nabla{\varphi}\mathrm{d}x=%
%TCIMACRO{\dint \nolimits_{B_{R}(x_{0})}}%
%BeginExpansion
{\displaystyle\int\nolimits_{B_{R}(x_{0})}}
%EndExpansion
\left(  f(w)+a(\left\vert x-x_{0}\right\vert )\left\vert \nabla w\right\vert
^{p-1}\right)  {\varphi}\mathrm{d}x.
\]
This shows that $w$ is also a weak solution to (\ref{Rf}).
\end{proof}

\subsection{Proofs of Theorems \ref{teovqa} and \ref{teowa}\label{teosvw}}

Let us fix $q\in\lbrack1,p)$ and search for a positive solution to the
boundary value problem%
\begin{equation}
\left\{
\begin{array}
[c]{l}%
\left(  e^{-A(r)}r^{N-1}\left\vert v^{\prime}\right\vert ^{p-1}\right)
^{\prime}=e^{-A(r)}r^{N-1}\left\vert v\right\vert ^{q-2}v,\text{ \ }%
r\in(0,R)\\
v^{\prime}(0)=0=v(R),
\end{array}
\right.  \label{erqq}%
\end{equation}
as a fixed point of the operator $T_{q}:C\left(  [0,R]\right)  \rightarrow
C\left(  [0,R]\right)  $, defined as%

\begin{equation}
\left(  T_{q}v\right)  (r)=\int_{r}^{R}\left(
%TCIMACRO{\dint \nolimits_{0}^{t}}%
%BeginExpansion
{\displaystyle\int\nolimits_{0}^{t}}
%EndExpansion
e^{A(t)-A(s)}\left(  \frac{s}{t}\right)  ^{N-1}\left\vert v(s)\right\vert
^{q-1}\mathrm{d}s\right)  ^{\frac{1}{p-1}}\mathrm{d}t,\text{ }0\leq r\leq R.
\label{Tq}%
\end{equation}

For each $0\leq r\leq R,$ we have the estimates%
\begin{align*}
\left(  T_{q}v\right)  (r)  &  \leq\beta^{\frac{1}{p-1}}\left\Vert
v\right\Vert _{\infty}^{\frac{q-1}{p-1}}\int_{r}^{R}\left(
%TCIMACRO{\dint \nolimits_{0}^{t}}%
%BeginExpansion
{\displaystyle\int\nolimits_{0}^{t}}
%EndExpansion
\left(  \frac{s}{t}\right)  ^{N-1}\mathrm{d}s\right)  ^{\frac{1}{p-1}%
}\mathrm{d}t\\
&  =\left(  \frac{\beta}{N}\right)  ^{\frac{1}{p-1}}\left\Vert v\right\Vert
_{\infty}^{\frac{q-1}{p-1}}\left(  \frac{p-1}{p}\right)  \left(  R^{\frac
{p}{p-1}}-r^{\frac{p}{p-1}}\right)
\end{align*}
and%
\[
\left\vert \left(  T_{q}v\right)  ^{\prime}(r)\right\vert \leq\left(
\frac{\beta}{r^{N-1}}%
%TCIMACRO{\dint \nolimits_{0}^{r}}%
%BeginExpansion
{\displaystyle\int\nolimits_{0}^{r}}
%EndExpansion
s^{N-1}\mathrm{d}s\right)  ^{\frac{1}{p-1}}\left\Vert v\right\Vert _{\infty
}^{\frac{q-1}{p-1}}=\left(  \frac{\beta r}{N}\right)  ^{\frac{1}{p-1}%
}\left\Vert v\right\Vert _{\infty}^{\frac{q-1}{p-1}}.
\]
It follows that
\begin{equation}
0\leq\left(  T_{q}v\right)  (r)\leq\left(  \frac{\beta R^{p}}{N}\right)
^{\frac{1}{p-1}}\left(  \frac{p-1}{p}\right)  \left\Vert v\right\Vert
_{\infty}^{\frac{q-1}{p-1}},\text{ }0\leq r\leq R, \label{ax1}%
\end{equation}
and%
\begin{equation}
\left\Vert \left(  T_{q}v\right)  ^{\prime}\right\Vert _{\infty}\leq\left(
\frac{\beta R}{N}\right)  ^{\frac{1}{p-1}}\left\Vert v\right\Vert _{\infty
}^{\frac{q-1}{p-1}}. \label{ax0}%
\end{equation}

It is well known that $\vartheta_{q}$ is radially symmetric and radially
decreasing. Therefore, $\vartheta_{q}(x)=\theta_{q}(\left\vert x-x_{0}%
\right\vert )$ where
\[
\left\{
\begin{array}
[c]{lll}%
\left(  r^{N-1}\left\vert \theta_{q}^{\prime}(r)\right\vert ^{p-1}\right)
^{\prime}=r^{N-1}\theta_{q}(r)^{q-1} & \mathrm{if} & 0<r<R\\
\theta_{q}^{\prime}(0)=0=\theta_{q}(R). &  &
\end{array}
\right.
\]

Thus, $\theta_{q}$ is a fixed point of $T_{q}$ in the particular case where
$a\equiv0,$ that is,%
\begin{equation}
\theta_{q}(r)=\int_{r}^{R}\left(
%TCIMACRO{\dint \nolimits_{0}^{t}}%
%BeginExpansion
{\displaystyle\int\nolimits_{0}^{t}}
%EndExpansion
\left(  \frac{s}{t}\right)  ^{N-1}\theta_{q}(s)^{q-1}\mathrm{d}s\right)
^{\frac{1}{p-1}}\mathrm{d}t. \label{rvq}%
\end{equation}
It follows that%
\[
\left\Vert \theta_{q}\right\Vert _{\infty}\leq\left\Vert \theta_{q}\right\Vert
_{\infty}^{\frac{q-1}{p-1}}\int_{r}^{R}\left(
%TCIMACRO{\dint \nolimits_{0}^{t}}%
%BeginExpansion
{\displaystyle\int\nolimits_{0}^{t}}
%EndExpansion
\left(  \frac{s}{t}\right)  ^{N-1}\mathrm{d}s\right)  ^{\frac{1}{p-1}%
}\mathrm{d}t\leq\left\Vert \theta_{q}\right\Vert _{\infty}^{\frac{q-1}{p-1}%
}\frac{p-1}{p}\left(  \frac{R^{p}}{N}\right)  ^{\frac{1}{p-1}},
\]
so that%
\[
\left\Vert \theta_{q}\right\Vert _{\infty}\leq\left(  \frac{p-1}{p}\right)
^{\frac{p-1}{p-q}}\left(  \frac{R^{p}}{N}\right)  ^{\frac{1}{p-q}}.
\]

Let us define the subset
\[
F:=\left\{  v\in C\left(  [0,R]\right)  :k\theta_{q}(r)\leq v(r)\leq M,\text{
for all }0\leq r\leq R.\right\}  ,
\]
where
\begin{equation}
k:=\alpha^{\frac{1}{p-q}}\text{ \ and \ }M:=\left(  \frac{p-1}{p}\right)
^{\frac{p-1}{p-q}}\left(  \frac{\beta R^{p}}{N}\right)  ^{\frac{1}{p-q}}.
\label{kM}%
\end{equation}

Taking into account that $0<\alpha<\beta,$ we have that the function
$k\theta_{q}\in F$ since%
\begin{align*}
k\theta_{q}(r)  &  \leq\alpha^{\frac{1}{p-q}}\left(  \frac{p-1}{p}\right)
^{\frac{p-1}{p-q}}\left(  \frac{R^{p}}{N}\right)  ^{\frac{1}{p-q}}\\
&  <\beta^{\frac{1}{p-q}}\left(  \frac{p-1}{p}\right)  ^{\frac{p-1}{p-q}%
}\left(  \frac{R^{p}}{N}\right)  ^{\frac{1}{p-q}}=M,
\end{align*}
for each $0\leq r\leq R.$

In order to show that $T_{q}(F)\subset F$ let us fix $v\in F.$ It follows from
(\ref{ax1}) that%
\[
\left(  T_{q}v\right)  (r)\leq\left(  \frac{\beta R^{p}}{N}\right)  ^{\frac
{1}{p-1}}\left(  \frac{p-1}{p}\right)  M^{\frac{q-1}{p-1}}=M,
\]
and from (\ref{rvq}) that
\begin{align*}
\left(  T_{q}v\right)  (r)  &  \geq\int_{r}^{R}\left(
%TCIMACRO{\dint \nolimits_{0}^{t}}%
%BeginExpansion
{\displaystyle\int\nolimits_{0}^{t}}
%EndExpansion
e^{A(t)-A(s)}\left(  \frac{s}{t}\right)  ^{N-1}\left(  k\theta_{q}(s)\right)
^{q-1}\mathrm{d}s\right)  ^{\frac{1}{p-1}}\mathrm{d}t\\
&  \geq\alpha^{\frac{1}{p-1}}k^{\frac{q-1}{p-1}}\int_{r}^{R}\left(
%TCIMACRO{\dint \nolimits_{0}^{t}}%
%BeginExpansion
{\displaystyle\int\nolimits_{0}^{t}}
%EndExpansion
\left(  \frac{s}{t}\right)  ^{N-1}\theta_{q}(s)^{q-1}\mathrm{d}s\right)
^{\frac{1}{p-1}}\mathrm{d}t\\
&  =\alpha^{\frac{1}{p-1}}k^{\frac{q-1}{p-1}}\theta_{q}(r)=k\theta_{q}(r).
\end{align*}

\begin{proposition}
\label{wn}The sequence $\left(  w_{n}\right)  \subset C([0,R]),$ defined
recursively as%
\[
w_{1}:=\alpha^{\frac{1}{p-q}}\theta_{q}\text{ \ and \ }w_{n+1}=T_{q}w_{n},
\]
is nondecreasing and converges uniformly to a function $\omega_{q}\in F\cap
C^{1}\left(  [0,R]\right)  \cap C^{2}\left(  (0,R]\right)  $ that solves
(\ref{erqq}). Moreover,
\begin{equation}
\alpha^{\frac{1}{p-q}}\theta_{q}(r)\leq\omega_{q}(r)\leq\left(  \frac{p-1}%
{p}\right)  ^{\frac{p-1}{p-q}}\left(  \frac{\beta R^{p}}{N}\right)  ^{\frac
{1}{p-q}}\text{ \ for all }r\in\lbrack0,R] \label{estw1}%
\end{equation}
and%
\begin{equation}
\left\Vert \omega_{q}^{\prime}\right\Vert _{\infty}\leq\left(  \frac{p-1}%
{p}\right)  ^{\frac{q-1}{p-q}}\left(  \frac{\beta R^{q}}{N}\right)  ^{\frac
{1}{p-q}}. \label{estw2}%
\end{equation}

\end{proposition}

\begin{proof}
Let us use the notation (\ref{kM}) and the fact that $T_{q}(F)\subset F.$ As
$w_{1}=k\theta_{q}\in F$ we have that $w_{2}=T_{q}w_{1}\in F,$ so that
\[
w_{1}(r)=k\theta_{q}\leq w_{2}(r)\leq M,\text{ \ }r\in\left[  0,R\right]  .
\]

Taking into account that $w_{n+1}=T_{q}w_{n}$ and recalling the definition of
$T_{q}$ in (\ref{Tq}), we argue by induction (using that $k\theta_{q}\leq
w_{n-1}\leq w_{n}\leq M$) to conclude that
\begin{equation}
k\theta_{q}\leq w_{n}\leq w_{n+1}\leq M\text{ \ in }\left[  0,R\right]  .
\label{estwn}%
\end{equation}

In addition, (\ref{ax0}) yields%
\begin{align*}
\left\Vert w_{n+1}^{\prime}\right\Vert _{\infty}  &  \leq\left(  \frac{\beta
R}{N}\right)  ^{\frac{1}{p-1}}\left\Vert w_{n}\right\Vert _{\infty}%
^{\frac{q-1}{p-1}}\\
&  \leq\left(  \frac{\beta R}{N}\right)  ^{\frac{1}{p-1}}M^{\frac{q-1}{p-1}%
}=\left(  \frac{p-1}{p}\right)  ^{\frac{q-1}{p-q}}\left(  \frac{\beta R^{q}%
}{N}\right)  ^{\frac{1}{p-q}}.
\end{align*}

Summarizing, the sequence $\left(  w_{n}\right)  $ is nondecreasing and
uniformly bounded in $C^{1}([0,R]).$ Consequently, $\left(  w_{n}\right)  $
converges uniformly in $\left[  0,R\right]  $ to a function $\omega_{q}\in
C\left(  \left[  0,R\right]  \right)  .$ Thus, we can to pass to the limit in
the identity%
\[
w_{n+1}(r)=\int_{r}^{R}\left(
%TCIMACRO{\dint \nolimits_{0}^{t}}%
%BeginExpansion
{\displaystyle\int\nolimits_{0}^{t}}
%EndExpansion
e^{A(t)-A(s)}\left(  \frac{s}{t}\right)  ^{N-1}w_{n}\left(  s\right)
^{q-1}\mathrm{d}s\right)  ^{\frac{1}{p-1}}\mathrm{d}t,
\]
and also in (\ref{estwn}) to conclude that $T\omega_{q}=\omega_{q}\in F.$
Therefore, $\omega_{q}$ is the positive solution to (\ref{erqq}), safisfies
(\ref{estw1}) and belongs to $C^{1}\left(  \left[  0,R\right]  \right)  .$
Furthermore, the estimates (\ref{ax0}) and (\ref{estw1}) yield%
\begin{align*}
\left\Vert \omega_{q}^{\prime}\right\Vert _{\infty}  &  \leq\left(
\frac{\beta R}{N}\right)  ^{\frac{1}{p-1}}\left\Vert \omega_{q}\right\Vert
_{\infty}^{\frac{q-1}{p-1}}\\
&  \leq\left(  \frac{\beta R}{N}\right)  ^{\frac{1}{p-1}}M^{\frac{q-1}{p-1}%
}\leq\left(  \frac{p-1}{p}\right)  ^{\frac{q-1}{p-q}}\left(  \frac{\beta
R^{q}}{N}\right)  ^{\frac{1}{p-q}}.
\end{align*}

\end{proof}

We now are ready to prove Theorems \ref{teovqa} and
\ref{teowa}.

\begin{proof}
[Proof of Theorem \ref{teovqa}]Let us define%
\[
w_{q}(x):=\omega_{q}\left(  \left\vert x-x_{0}\right\vert \right)  ,\text{
\ }x\in B_{R}(x_{0}),
\]
where $\omega_{q}$ is the solution to (\ref{erqq}) given by Proposition
\ref{wn}. It follows that $w_{q}\in C^{2}\left(  B_{R}(x_{0})\setminus\left\{
x_{0}\right\}  \right)  \cap C_{0}^{1}\left(  \overline{B_{R}(x_{0})}\right)
$ is radially symmetric and radially decreasing. According to Lemma \ref{qw}
applyied to $f(t)=\left\vert t\right\vert ^{q-2}t,$ the function $w_{q}$ is a
positive weak solution to both problems (\ref{eRq}) and (\ref{Rq}). By
uniqueness of the positive solutions to (\ref{eRq}), we conclude that
$w_{q}=w_{q,a}.$ This proves that $w_{q,a}$ is radially symmetric and radially
decreasing. The bounds (\ref{estvqa1}) and (\ref{estvqa2}) follow respectively
from (\ref{estw1}) and (\ref{estw2}).

We remark that $w_{q,a}$ is the only positive, radially symmetric weak
solution to (\ref{Rq}) since any positive weak solution to (\ref{Rq}) with
this symmetry is also a positive weak solution to (\ref{eRq}).
\end{proof}

\begin{proof}
[Proof of Theorem \ref{teowa}]Let us define
\[
\upsilon_{q}(r):=\frac{\omega_{q}(r)}{\left\Vert \omega_{q}\right\Vert
_{\infty}}\text{ \ and \ }\gamma_{q}:=\frac{1}{\left\Vert \omega
_{q}\right\Vert _{\infty}^{p-q}},
\]
where $\omega_{q}$ is the solution to (\ref{erqq}) given by Proposition
\ref{wn}. It follows that $\left\Vert \upsilon_{q}\right\Vert _{\infty}=1,$%
\begin{equation}
\upsilon_{q}(r)=\int_{r}^{R}\left(
%TCIMACRO{\dint \nolimits_{0}^{t}}%
%BeginExpansion
{\displaystyle\int\nolimits_{0}^{t}}
%EndExpansion
e^{A(t)-A(s)}\left(  \frac{s}{t}\right)  ^{N-1}\gamma_{q}\upsilon_{q}\left(
s\right)  ^{q-1}\mathrm{d}s\right)  ^{\frac{1}{p-1}}\mathrm{d}t\label{vqp}%
\end{equation}
and%
\[
\left(  \frac{p}{p-1}\right)  ^{p-1}\frac{N}{\beta R^{p}}\leq\gamma_{q}%
\leq\frac{1}{\alpha\left\Vert \vartheta_{q}\right\Vert _{\infty}^{q-p}}%
\]
(recall that $\vartheta_{q}(x)=\theta_{q}(\left\vert x-x_{0}\right\vert )$).

We recall from Theorem \ref{sublin} that
\[
\lim_{q\rightarrow p^{-}}\frac{1}{\left\Vert \vartheta_{q}\right\Vert
_{\infty}^{q-p}}=\lambda_{1,p}\left(  B_{R}(x_{0})\right)  .
\]
Hence,
\[
\left(  \frac{p}{p-1}\right)  ^{p-1}\frac{N}{\beta R^{p}}\leq\liminf
_{q\rightarrow p^{-}}\gamma_{q}\leq\limsup_{q\rightarrow p^{-}}\gamma_{q}%
\leq\frac{\lambda_{1,p}\left(  B_{R}(x_{0})\right)  }{\alpha}%
\]
and
\[
\limsup_{q\rightarrow p^{-}}\left\Vert \upsilon_{q}^{\prime}\right\Vert
_{\infty}\leq\limsup_{q\rightarrow p^{-}}\left(  \beta\gamma_{q}\frac{R}%
{N}\right)  ^{\frac{1}{p-1}}\leq\left(  \frac{\beta R}{\alpha N}\lambda
_{1,p}\left(  B_{R}(x_{0})\right)  \right)  ^{\frac{1}{p-1}},
\]
since%
\[
\left\vert \upsilon_{q}^{\prime}(r)\right\vert ^{p-1}=\gamma_{q}%
%TCIMACRO{\dint \nolimits_{0}^{r}}%
%BeginExpansion
{\displaystyle\int\nolimits_{0}^{r}}
%EndExpansion
e^{A(t)-A(s)}\left(  \frac{s}{r}\right)  ^{N-1}\upsilon_{q}\left(  s\right)
^{q-1}\mathrm{d}s\leq\beta\gamma_{q}\frac{R}{N}.
\]

The uniform boundedness of $\gamma_{q},$ $\left\Vert \upsilon_{q}\right\Vert
_{\infty}$ and $\left\Vert \upsilon_{q}^{\prime}\right\Vert _{\infty}$
guarantees that, up to a subsequence,
\[
\lim_{q\rightarrow p^{-}}\upsilon_{q}=\upsilon_{p},\text{ uniformly in
}\overline{B_{R}(x_0)}\text{,}%
\]
for some $\upsilon_{p}\in C\left(  \left[  0,R\right]  \right)  ,$ and
\[
\lim_{q\rightarrow p^{-}}\gamma_{q}=\gamma_{p},
\]
for some constant $\gamma_{p}$ satisfying
\begin{equation}
\left(  \frac{p}{p-1}\right)  ^{p-1}\frac{N}{\beta R^{p}}\leq\gamma_{p}%
\leq\frac{\lambda_{1,p}\left(  B_{R}(x_{0})\right)  }{\alpha}.\label{estgp}%
\end{equation}

Hence, after letting $q\rightarrow p^{-}$ in (\ref{vqp}), we find%
\begin{equation}
\upsilon_{p}(r)=\int_{r}^{R}\left(
%TCIMACRO{\dint \nolimits_{0}^{t}}%
%BeginExpansion
{\displaystyle\int\nolimits_{0}^{t}}
%EndExpansion
e^{A(t)-A(s)}\left(  \frac{s}{t}\right)  ^{N-1}\gamma_{p}\upsilon_{p}\left(
s\right)  ^{p-1}\mathrm{d}s\right)  ^{\frac{1}{p-1}}\mathrm{d}t.\label{vpr}%
\end{equation}
It follows that $\upsilon_{p}\in C^{1}\left(  \left[  0,R\right]  \right)
\cap C^{2}\left(  (0,R]\right)  $ is a positive solution to the ODE
boundary value problem%
\[
\left\{
\begin{array}
[c]{lll}%
\left(  e^{-A(r)}r^{N-1}\left\vert \upsilon^{\prime}\right\vert ^{p-1}\right)
^{\prime}=\gamma_{p}e^{-A(r)}r^{N-1}\left\vert \upsilon\right\vert
^{p-2}\upsilon & \mathrm{if} & r\in(0,R)\\
\upsilon^{\prime}(0)=0=\upsilon(R). &  &
\end{array}
\right.
\]
In view of Lemma \ref{qw} applyied to $f(t)=\left\vert t\right\vert ^{p-2}t,$
the function $\upsilon_{p}\left(  \left\vert x-x_{0}\right\vert \right)  $ is
a weak solution to both problems%
\begin{equation}
\left\{
\begin{array}
[c]{lll}%
-\Delta_{p}w=\gamma_{p}\left\vert w\right\vert ^{p-2}w+a(\left\vert
x-x_{0}\right\vert )\left\vert \nabla w\right\vert ^{p-1} & \mathrm{in} &
B_{R}(x_{0})\\
w=0 & \mathrm{on} & \partial B_{R}(x_{0})
\end{array}
\right.  \label{wax}%
\end{equation}
and%
\begin{equation}
\left\{
\begin{array}
[c]{lll}%
-\operatorname{div}\left(  e^{-A(\left\vert x-x_{0}\right\vert )}\left\vert
\nabla w\right\vert ^{p-2}\nabla w\right)  =\gamma_{p}e^{-A(\left\vert
x-x_{0}\right\vert )}\left\vert w\right\vert ^{p-2}w & \mathrm{in} &
B_{R}(x_{0})\\
w=0 & \mathrm{on} & \partial B_{R}(x_{0}).
\end{array}
\right.  \label{ewax}%
\end{equation}
According to Lemma \ref{lemwpa}, $\upsilon_{p}\left(  \left\vert
x-x_{0}\right\vert \right)  =w_{a}\left(  x\right)  $ and $\gamma_{p}%
=\sigma_{a}\left(  B_{R}(x_{0})\right)  .$ Hence, (\ref{estgp}) leads to
(\ref{estsa}). The estimate (\ref{estwa}) follows after deriving (\ref{vpr})
since%
\[
\left\Vert \nabla w_{a}\right\Vert _{\infty}^{p-1}=\max_{0\leq r\leq
R}\left\vert \upsilon_{p}^{\prime}(r)\right\vert ^{p-1}\leq\frac{\beta R}%
{N}\gamma_{p}=\frac{\beta R}{N}\sigma_{a}\left(  B_{R}(x_{0})\right)  .
\]

If $\widetilde{w}$ is a positive, radially symmetric weak solution to
(\ref{PBRl}), then $\widetilde{w}$ is a positive weak solution to the problem
\begin{equation}
\left\{
\begin{array}
[c]{lll}%
-\operatorname{div}\left(  e^{-A(\left\vert x-x_{0}\right\vert )}\left\vert
\nabla w\right\vert ^{p-2}\nabla w\right)  =\lambda e^{-A(\left\vert
x-x_{0}\right\vert )}\left\vert w\right\vert ^{p-2}w & \mathrm{in} &
B_{R}(x_{0})\\
w=0 & \mathrm{on} & \partial B_{R}(x_{0}).
\end{array}
\right.  \label{elwax}%
\end{equation}
Hence, Lemma \ref{lemwpa} again implies that $\lambda=\sigma_{a}(B_{R}%
(x_{0}))$ and $\widetilde{w}=cw_{a}$ for some constant $c>0.$
\end{proof}

\section{Proof of Theorem \ref{teomain}\label{secfinal}}

Let $\mu\in\mathbb{R}$ be a constant, $x_{0}\in\mathbb{R}^{N}$ and $R>0$. Let
us consider the problem%
\begin{equation}
\left\{
\begin{array}
[c]{lll}%
-\Delta_{p}u=\left\vert u\right\vert ^{q-2}u+\mu\left\vert \nabla u\right\vert
^{p-1} & \mathrm{in} & B_{R}(x_{0})\\
u=0 & \mathrm{on} & \partial B_{R}(x_{0}),
\end{array}
\right.  \label{Pmu}%
\end{equation}
which is (\ref{Rq}) with $a\equiv\mu.$ According to Theorem
\ref{teovqa}, problem (\ref{Pmu}) has a unique positive, radially symmetric
weak solution $u_{q,\mu}\in C^{1}\left(  \overline{B_{R}(x_{0})}\right)  \cap
C^{2}\left(  B_{R}(x_{0})\setminus\left\{  x_{0}\right\}  \right)  .$
Moreover,%
\begin{equation}
\alpha^{\frac{1}{p-q}}\vartheta_{q}\leq u_{q,\mu}\leq\left(  \frac{p-1}%
{p}\right)  ^{\frac{p-1}{p-q}}\left(  \frac{\beta R^{p}}{N}\right)  ^{\frac
{1}{p-q}} \label{uqmu}%
\end{equation}
and
\[
\left\Vert \nabla u_{q,\mu}\right\Vert _{\infty}\leq\left(  \frac{p-1}%
{p}\right)  ^{\frac{p-1}{p-q}}\left(  \frac{\beta R^{q}}{N}\right)  ^{\frac
{1}{p-q}},
\]
where
\begin{equation}
\alpha=\min\left\{  e^{\mu R},1\right\}  \text{ \ and \ }\beta=\max\left\{
e^{\mu R},1\right\}  . \label{ab1}%
\end{equation}

We recall that $u_{q,\mu}(x)  =\omega\left(  \left\vert
x-x_{0}\right\vert \right)  $ for all $x\in B_{R}(x_{0}),$ where $\omega\in
C^{2}\left(  (0,R]\right)  \cap C^{1}\left(  \left[  0,R\right]  \right)  .$
The properties of $\omega$ imply that%

\begin{equation}
-\operatorname{div}\left(  \left\vert \nabla u_{q,\mu}\right\vert ^{p-2}\nabla
u_{q,\mu}\right)  =\left\vert u_{q,\mu}\right\vert ^{q-2}u_{q,\mu}%
+\mu\left\vert \nabla u_{q,\mu}\right\vert ^{p-1}\text{ \ pointwise in }B_{R}%
(x_{0})\setminus\left\{  x_{0}\right\}  \label{-div}%
\end{equation}
and%
\begin{equation}
\nabla u_{q,\mu}(x)  =-\left\vert \omega^{\prime}(
\left\vert x-x_{0}\right\vert)  \right\vert \frac{x-x_{0}}{\left\vert
x-x_{0}\right\vert }\text{ \ for all }x\in B_{R}(x_{0})\setminus\left\{
x_{0}\right\}  . \label{gradu}%
\end{equation}

\begin{lemma}
\label{sub}Let $x_{0}\in\Omega$ and $\mu\in\mathbb{R}.$ Let $u_{q,\mu}$ be the
only positive, radially symmetric weak solution to the problem (\ref{Pmu}),
where $R>0$ is such that $B_{R}(x_{0})\subset\Omega.$ Define $\underline{u}\in
W_{0}^{1,p}(\Omega)$ as the extension of $u_{q,\mu}$ by zero to $\overline
{\Omega}\setminus B_{R}(x_{0}).$ Then $\underline{u}$ is a subsolution to the
problem%
\begin{equation}
\left\{
\begin{array}
[c]{lll}%
-\Delta_{p}u=\left\vert u\right\vert ^{q-2}u+\mu\left\vert \nabla u\right\vert
^{p-1} & \mathrm{in} & \Omega\\
u=0 & \mathrm{on} & \partial\Omega.
\end{array}
\right.  \label{sub2}%
\end{equation}

Moreover, if $\mu\leq0,$ then $\underline{u}\leq\vartheta_{q}$ in $\Omega$.
\end{lemma}

\begin{proof}
Let $\phi\in C_{c}^{\infty}(\Omega)$ be such that $\phi\geq0$ and take
$0<\epsilon<R.$ It follows that%
\begin{align*}%
%TCIMACRO{\dint \nolimits_{\Omega}}%
%BeginExpansion
{\displaystyle\int\nolimits_{\Omega}}
%EndExpansion
\left\vert \nabla\underline{u}\right\vert ^{p-2}\nabla\underline{u}\cdot
\nabla\phi\mathrm{d}x  &  =%
%TCIMACRO{\dint \nolimits_{B_{R}(x_{0})}}%
%BeginExpansion
{\displaystyle\int\nolimits_{B_{R}(x_{0})}}
%EndExpansion
\left\vert \nabla u_{q,\mu}\right\vert ^{p-2}\nabla u_{q,\mu}\cdot\nabla
\phi\mathrm{d}x\\
&  =I_{\epsilon}+K_{\epsilon}-L_{\epsilon}%
\end{align*}
where%
\[
I_{\epsilon}:=%
%TCIMACRO{\dint \nolimits_{B_{\epsilon}(x_{0})}}%
%BeginExpansion
{\displaystyle\int\nolimits_{B_{\epsilon}(x_{0})}}
%EndExpansion
\left\vert \nabla u_{q,\mu}\right\vert ^{p-2}\nabla u_{q,\mu}\cdot\nabla
\phi\mathrm{d}x,
\]%
\[
K_{\epsilon}:=%
%TCIMACRO{\dint \nolimits_{B_{R}(x_{0})\setminus B_{\epsilon}(x_{0})}}%
%BeginExpansion
{\displaystyle\int\nolimits_{B_{R}(x_{0})\setminus B_{\epsilon}(x_{0})}}
%EndExpansion
\operatorname{div}\left(  \phi\left\vert \nabla u_{q,\mu}\right\vert
^{p-2}\nabla u_{q,\mu}\right)  \mathrm{d}x
\]
and%
\[
L_{\epsilon}:=%
%TCIMACRO{\dint \nolimits_{B_{R}(x_{0})\setminus B_{\epsilon}(x_{0})}}%
%BeginExpansion
{\displaystyle\int\nolimits_{B_{R}(x_{0})\setminus B_{\epsilon}(x_{0})}}
%EndExpansion
\phi\operatorname{div}\left(  \left\vert \nabla u_{q,\mu}\right\vert
^{p-2}\nabla u_{q,\mu}\right)  \mathrm{d}x.
\]

Taking into account that
\[
\left\vert I_{\epsilon}\right\vert \leq%
%TCIMACRO{\dint \nolimits_{B_{\epsilon}(x_{0})}}%
%BeginExpansion
{\displaystyle\int\nolimits_{B_{\epsilon}(x_{0})}}
%EndExpansion
\left\vert \nabla u_{q,\mu}\right\vert ^{p-1}\left\vert \nabla\phi\right\vert
\mathrm{d}x\leq\left\Vert \nabla u_{q,\mu}\right\Vert _{\infty}^{p-1}%
%TCIMACRO{\dint \nolimits_{B_{\epsilon}(x_{0})}}%
%BeginExpansion
{\displaystyle\int\nolimits_{B_{\epsilon}(x_{0})}}
%EndExpansion
\left\vert \nabla\phi\right\vert \mathrm{d}x
\]
we conclude that $\lim\limits_{\epsilon\rightarrow0^{+}}I_{\epsilon}=0.$

Note from (\ref{gradu}) that
\[%
%TCIMACRO{\dint \nolimits_{\partial B_{r}(x_{0})}}%
%BeginExpansion
{\displaystyle\int\nolimits_{\partial B_{r}(x_{0})}}
%EndExpansion
\phi\left\vert \nabla u_{q,\mu}\right\vert ^{p-2}\nabla u_{q,\mu}\cdot
\frac{x-x_{0}}{\left\vert x-x_{0}\right\vert }\mathrm{d}S_{x}=-\left\vert
\omega^{\prime}(r)\right\vert ^{p-1}%
%TCIMACRO{\dint \nolimits_{\partial B_{r}(x_{0})}}%
%BeginExpansion
{\displaystyle\int\nolimits_{\partial B_{r}(x_{0})}}
%EndExpansion
\phi\left(  x\right)  \mathrm{d}S_{x}%
\]
for all $\epsilon<r<R.$ Hence, owing to the Divergence Theorem,%
\[
K_{\epsilon}=-\left\vert \omega^{\prime}(R)  \right\vert ^{p-1}%
%TCIMACRO{\dint \nolimits_{\partial B_{R}(x_{0})}}%
%BeginExpansion
{\displaystyle\int\nolimits_{\partial B_{R}(x_{0})}}
%EndExpansion
\phi(x)  \mathrm{d}S_{x}+\left\vert \omega^{\prime}(
\epsilon)  \right\vert ^{p-1}%
%TCIMACRO{\dint \nolimits_{\partial B_{\epsilon}(x_{0})}}%
%BeginExpansion
{\displaystyle\int\nolimits_{\partial B_{\epsilon}(x_{0})}}
%EndExpansion
\phi\mathrm{d}S_{x},%
\]
so that
\[
\lim_{\epsilon\rightarrow0^{+}}K_{\epsilon}=-\left\vert \omega^{\prime}(
R)  \right\vert ^{p-1}%
%TCIMACRO{\dint \nolimits_{\partial B_{R}(x_{0})}}%
%BeginExpansion
{\displaystyle\int\nolimits_{\partial B_{R}(x_{0})}}
%EndExpansion
\phi(x)  \mathrm{d}S_{x}\leq 0,
\]
since
\[
\lim_{\epsilon\rightarrow0^{+}}\left\vert \omega^{\prime}(
\epsilon)  \right\vert ^{p-1}%
%TCIMACRO{\dint \nolimits_{\partial B_{\epsilon}(x_{0})}}%
%BeginExpansion
{\displaystyle\int\nolimits_{\partial B_{\epsilon}(x_{0})}}
%EndExpansion
\phi (x)\mathrm{d}S_{x}=0.
\]

It follows from (\ref{-div}) that%
\[
L_{\epsilon}=-%
%TCIMACRO{\dint \nolimits_{B_{R}(x_{0})\setminus B_{\epsilon}(x_{0})}}%
%BeginExpansion
{\displaystyle\int\nolimits_{B_{R}(x_{0})\setminus B_{\epsilon}(x_{0})}}
%EndExpansion
\left[  \left\vert u_{q,\mu}\right\vert ^{q-2}u_{q,\mu}+\mu\left\vert \nabla
u_{q,\mu}\right\vert ^{p-1}\right]  \phi\mathrm{d}x.
\]
Hence, as $\left\vert u_{q,\mu}\right\vert ^{q-2}u_{q,\mu}+\mu\left\vert
\nabla u_{q,\mu}\right\vert ^{p-1}\in C\left(  \overline{B_{R}(x_{0})}\right)
,$ we have that
\begin{align*}
\lim_{\epsilon\rightarrow0^{+}}L_{\epsilon}  &  =-%
%TCIMACRO{\dint \nolimits_{B_{R}(x_{0})}}%
%BeginExpansion
{\displaystyle\int\nolimits_{B_{R}(x_{0})}}
%EndExpansion
\left[  \left\vert u_{q,\mu}\right\vert ^{q-2}u_{q,\mu}+\mu\left\vert \nabla
u_{q,\mu}\right\vert ^{p-1}\right]  \phi\mathrm{d}x\\
&  =-%
%TCIMACRO{\dint \nolimits_{\Omega}}%
%BeginExpansion
{\displaystyle\int\nolimits_{\Omega}}
%EndExpansion
\left[  \left\vert \underline{u}\right\vert ^{q-2}\underline{u}+\mu\left\vert
\nabla\underline{u}\right\vert ^{p-1}\right]  \phi\mathrm{d}x.
\end{align*}

Therefore,
\begin{align*}%
%TCIMACRO{\dint \nolimits_{\Omega}}%
%BeginExpansion
{\displaystyle\int\nolimits_{\Omega}}
%EndExpansion
\left\vert \nabla\underline{u}\right\vert ^{p-2}\nabla\underline{u}\cdot
\nabla\phi\mathrm{d}x  &  =\lim_{\epsilon\rightarrow0^{+}}\left(  I_{\epsilon
}+K_{\epsilon}-L_{\epsilon}\right) \\
&  =-\left\vert \omega^{\prime}(R)  \right\vert ^{p-1}%
%TCIMACRO{\dint \nolimits_{\partial B_{R}(x_{0})}}%
%BeginExpansion
{\displaystyle\int\nolimits_{\partial B_{R}(x_{0})}}
%EndExpansion
\phi\left(  x\right)  \mathrm{d}S_{x}+%
%TCIMACRO{\dint \nolimits_{\Omega}}%
%BeginExpansion
{\displaystyle\int\nolimits_{\Omega}}
%EndExpansion
\left[  \left\vert \underline{u}\right\vert ^{q-2}\underline{u}+\mu\left\vert
\nabla\underline{u}\right\vert ^{p-1}\right]  \phi\mathrm{d}x\\
&  \leq%
%TCIMACRO{\dint \nolimits_{\Omega}}%
%BeginExpansion
{\displaystyle\int\nolimits_{\Omega}}
%EndExpansion
\left[  \left\vert \underline{u}\right\vert ^{q-2}\underline{u}+\mu\left\vert
\nabla\underline{u}\right\vert ^{p-1}\right]  \phi\mathrm{d}x,
\end{align*}
which, combined with the fact that $\underline{u}=0$ on $\partial\Omega,$
implies that $\underline{u}$ is a subsolution to (\ref{sub2}).

In particular, if $\mu\leq0,$ then $-\Delta_{p}\underline{u}\leq\underline
{u}^{q-1}$ in $\Omega,$ since%
\[%
%TCIMACRO{\dint \nolimits_{\Omega}}%
%BeginExpansion
{\displaystyle\int\nolimits_{\Omega}}
%EndExpansion
\left\vert \nabla\underline{u}\right\vert ^{p-2}\nabla\underline{u}\cdot
\nabla\phi\mathrm{d}x\leq%
%TCIMACRO{\dint \nolimits_{\Omega}}%
%BeginExpansion
{\displaystyle\int\nolimits_{\Omega}}
%EndExpansion
\left[  \left\vert \underline{u}\right\vert ^{q-2}\underline{u}+\mu\left\vert
\nabla\underline{u}\right\vert ^{p-1}\right]  \phi\mathrm{d}x\leq%
%TCIMACRO{\dint \nolimits_{\Omega}}%
%BeginExpansion
{\displaystyle\int\nolimits_{\Omega}}
%EndExpansion
\left\vert \underline{u}\right\vert ^{q-2}\underline{u}\phi\mathrm{d}x,
\]
for all $\phi\in W_{0}^{1,p}(\Omega),$ $\phi\geq0.$ As\ $\underline{u}\in
W^{1,p}(\Omega)\cap C(\overline{\Omega}),$ Corollary \ref{abpe2} implies that
$\underline{u}\leq\vartheta_{q}$ in $\Omega.$
\end{proof}

\begin{proof}
[Proof of Theorem \ref{teomain}]Let $x_{0}\in\Omega,$ $x_{1}\in\mathbb{R}^{N}$
be such that
\[
B_{r_{\Omega}}(x_{0})\subset\Omega\subset B_{R_{\Omega}}(x_{1}).
\]

Let us set
\[
\mu_{0}:=-\left\Vert a^{-}\right\Vert _{\infty}\text{ \ and \ }\mu
_{1}:=\left\Vert a^{+}\right\Vert _{\infty}.
\]
Note that
\begin{equation}
\mu_{0}\leq-a^{-}(x)\leq a(x)\leq a^{+}(x)\leq\mu_{1},
\label{ineq}%
\end{equation}
for all $x\in\Omega.$

Let $u_{q,\mu_{0}}\in W_{0}^{1,p}(B_{r_{\Omega}}(x_{0}))$ be the positive,
radially symmetric weak solution to the problem%
\[
\left\{
\begin{array}
[c]{lll}%
-\Delta_{p}u=\left\vert u\right\vert ^{q-2}u+\mu_{0}\left\vert \nabla
u\right\vert ^{p-1} & \mathrm{in} & B_{r_{\Omega}}(x_{0})\\
u=0 & \mathrm{on} & \partial B_{r_{\Omega}}(x_{0})
\end{array}
\right.
\]
and let $u_{q,\mu_{1}}\in W_{0}^{1,p}(B_{R_{\Omega}}(x_{1}))$ be the positive,
radially symmetric weak solution to the problem%
\[
\left\{
\begin{array}
[c]{lll}%
-\Delta_{p}u=\left\vert u\right\vert ^{q-2}u+\mu_{1}\left\vert \nabla
u\right\vert ^{p-1} & \mathrm{in} & B_{R_{\Omega}}(x_{1})\\
u=0 & \mathrm{on} & \partial B_{R_{\Omega}}(x_{1}).
\end{array}
\right.
\]

According to (\ref{uqmu}) and (\ref{ab1})
\[
u_{q,\mu_{1}}(x)  \leq\left(  \frac{p-1}{p}\right)  ^{\frac
{p-1}{p-q}}\left(  \frac{R_{\Omega}^{p}e^{\mu_{1}R_{\Omega}}}{N}\right)
^{\frac{1}{p-q}}\text{ \ for all }x\in B_{R_{\Omega}}(x_{1})
\]
and
\[
e^{\frac{-\left\vert \mu_{0}\right\vert r_{\Omega}}{p-q}}\vartheta
_{0,q}(x)  \leq u_{q,\mu_{0}}(x)  \text{ \ for all
}x\in B_{r_{\Omega}}(x_{0}),
\]
where $\vartheta_{0,q}$ denotes the only positive solution to the problem
(\ref{vq}) in $B_{r_{\Omega}}(x_{0}).$ Note that
\[
e^{\mu_{1}R_{\Omega}}=\max\left\{  e^{\mu_{1}R_{\Omega}},1\right\}  \text{
\ and \ }e^{-\left\vert \mu_{0}\right\vert r_{\Omega}}=\min\left\{  e^{\mu
_{0}r_{\Omega}},1\right\}  .
\]

Let $\underline{u}_{q}\in W_{0}^{1,p}(\Omega)$ be the extension of
$u_{q,\mu_{0}}$ by zero in $\overline{\Omega}\setminus B_{r_{\Omega}}(x_{0}).$
It follows from Lemma \ref{sub} that $\underline{u}_{q}$ is a subsolution to
the problem
\[
\left\{
\begin{array}
[c]{lll}%
-\Delta_{p}u=\left\vert u\right\vert ^{q-2}u+\mu_{0}\left\vert \nabla
u\right\vert ^{p-1} & \mathrm{in} & \Omega\\
u=0 & \mathrm{on} & \partial\Omega
\end{array}
\right.
\]
and $\underline{u}_{q}\leq\vartheta_{q}$ \ in $\Omega.$
Owing to (\ref{ineq}) we have that
\[
\left\vert \underline{u}_{q}\right\vert ^{q-2}\underline{u}_{q}+\mu
_{0}\left\vert \nabla\underline{u}_{q}\right\vert ^{p-1}\leq\left\vert
\underline{u}_{q}\right\vert ^{q-2}\underline{u}_{q}+a(x)\left\vert
\nabla\underline{u}_{q}\right\vert ^{p-1}\text{ in }\Omega.
\]
This shows that $\underline{u}_{q}$ is also a nonnegative subsolution to
problem
\begin{equation}
\left\{
\begin{array}
[c]{lll}%
-\Delta_{p}w=\left\vert w\right\vert ^{q-2}w+a(x)\left\vert \nabla
w\right\vert ^{p-1} & \mathrm{in} & \Omega\\
w=0 & \mathrm{on} & \partial\Omega.
\end{array}
\right.  \label{aq}%
\end{equation}

The function $u_{q,\mu_{1}}$ restricted to $\Omega$ is naturally a positive
supersolution to the problem%
\[
\left\{
\begin{array}
[c]{lll}%
-\Delta_{p}u=\left\vert u\right\vert ^{q-2}u+\mu_{1}\left\vert \nabla
u\right\vert ^{p-1} & \mathrm{in} & \Omega\\
u=0 & \mathrm{on} & \partial\Omega.
\end{array}
\right.
\]
Moreover, $-\Delta_{p}u_{q,\mu_{1}
}\geq u_{q,\mu_{1}}^{q-1}$ in $\Omega$ since $\mu_{1}\geq0$. Hence, as $u_{q,\mu_{1}}\in
W^{1,p}(\Omega)\cap C^{1}(\Omega)$ and $u_{q,\mu_{1}}>0$ in $\Omega,$ Theorem
\ref{Ape} (with $f(u)=u^{q-1}$) yields $\vartheta_{q}\leq u_{q,\mu_{1}}$\ in
$\Omega.$

It follows that the functions $\underline{u}_{q}\ $and $u_{q,\mu_{1}}$ are
ordered:%
\[
\underline{u}_{q}\leq\vartheta_{q}\leq u_{q,\mu_{1}}\text{ \ in }\Omega.
\]

According to (\ref{ineq}) we have that
\[
\left\vert u_{q,\mu_{1}}\right\vert ^{q-2}u_{q,\mu_{1}}+\mu_{1}\left\vert
\nabla u_{q,\mu_{1}}\right\vert ^{p-1}\geq\left\vert u_{q,\mu_{1}}\right\vert
^{q-2}u_{q,\mu_{1}}+a(x)\left\vert \nabla u_{q,\mu_{1}}\right\vert
^{p-1}\text{ in }\Omega.
\]
This implies that $u_{q,\mu_{1}}$ restricted to $\Omega$ is a supersolution
to the problem (\ref{aq}).

Summarizing, we have shown that $\underline{u}_{q}$ and $u_{q,\mu_{1}}$ form
an ordered pair of sub-supersolutions to (\ref{aq}). Moreover,
\begin{align*}
e^{\frac{-\left\vert \mu_{0}\right\vert r_{\Omega}}{p-q}}\underline{\vartheta
}_{0,q}(x) &  \leq\underline{u}_{q}(x)\\
&  \leq\vartheta_{q}(x)\leq u_{q,\mu_{1}}(x)\leq\left(  \frac{p-1}{p}\right)
^{\frac{p-1}{p-q}}\left(  \frac{R_{\Omega}^{p}e^{\mu_{1}R_{\Omega}}}%
{N}\right)  ^{\frac{1}{p-q}}\text{ \ for all }x\in\Omega,
\end{align*}
where $\underline{\vartheta}_{0,q}$ denotes the extension by zero of
$\vartheta_{0,q}.$ These bounds imply that%
\begin{equation}
e^{-\left\vert \mu_{0}\right\vert r_{\Omega}}\left\Vert \vartheta
_{0,q}\right\Vert _{\infty}^{p-q}\leq\left\Vert \underline{u}_{q}\right\Vert
_{\infty}^{p-q}\leq\left\Vert u_{q,\mu_{1}}\right\Vert _{\infty}^{p-q}%
\leq\left(  \frac{p-1}{p}\right)  ^{p-1}\frac{R_{\Omega}^{p}e^{\mu
_{1}R_{\Omega}}}{N}\label{bgam}%
\end{equation}
(note that $\left\Vert \underline{\vartheta}_{0,q}\right\Vert _{\infty
}=\left\Vert \vartheta_{0,q}\right\Vert _{\infty}$).

Therefore, according to Theorem \ref{Leon}, problem (\ref{aq}) has a weak
solution $w_{q}\in W_{0}^{1,p}(\Omega)$ between $\underline{u}_{q}$ and
$u_{q,\mu_{1}}.$

Let us define%
\[
u_{q}:=\frac{w_{q}}{\left\Vert w_{q}\right\Vert _{\infty}}\text{ \ and
\ }\gamma_{q}:=\frac{1}{\left\Vert w_{q}\right\Vert _{\infty}^{q-p}}.
\]
It follows that%
\[
\left\{
\begin{array}
[c]{lll}%
-\Delta_{p}u_{q}=\gamma_{q}\left\vert u_{q}\right\vert ^{q-2}u_{q}%
+a(x)\left\vert \nabla u_{q}\right\vert ^{p-1} & \mathrm{in} & \Omega\\
u_{q}=0 & \mathrm{on} & \partial\Omega
\end{array}
\right.
\]
in the weak sense, that is,%
\begin{equation}%
%TCIMACRO{\dint \nolimits_{\Omega}}%
%BeginExpansion
{\displaystyle\int\nolimits_{\Omega}}
%EndExpansion
\left\vert \nabla u_{q}\right\vert ^{p-2}\nabla u_{q}\cdot\nabla\phi
\mathrm{d}x=%
%TCIMACRO{\dint \nolimits_{\Omega}}%
%BeginExpansion
{\displaystyle\int\nolimits_{\Omega}}
%EndExpansion
\left[  \gamma_{q}\left\vert u_{q}\right\vert ^{q-2}u_{q}+a(x)\left\vert
\nabla u_{q}\right\vert ^{p-1}\right]  \phi\mathrm{d}x \label{weak1}%
\end{equation}
for all $\phi\in W_{0}^{1,p}(\Omega).$

The bounds (\ref{bgam}) imply that%
\[
\left(  \frac{p}{p-1}\right)  ^{p-1}\frac{Ne^{-\mu_{1}R_{\Omega}}}{R_{\Omega
}^{p}}\leq\gamma_{q}\leq\frac{e^{\left\vert \mu_{0}\right\vert r_{\Omega}}%
}{\left\Vert \vartheta_{0,q}\right\Vert _{\infty}^{p-q}}.
\]
Thus, taking into account that
\[
\lim_{q\rightarrow p^{-}}\frac{1}{\left\Vert \vartheta_{0,q}\right\Vert
_{\infty}^{p-q}}=\lambda_{1,p}\left(  B_{r_{\Omega}}(x_{0})\right)
\]
we have that, up to a subsequence,
\[
\lim_{q\rightarrow p^{-}}\gamma_{q}=\lambda_{a},
\]
for some $\lambda_{a}$ satisfying%
\begin{equation}
\left(  \frac{p}{p-1}\right)  ^{p-1}\frac{Ne^{-\mu_{1}R_{\Omega}}}{R_{\Omega
}^{p}}\leq\lambda_{a}\leq e^{\left\vert \mu_{0}\right\vert r_{\Omega}}%
\lambda_{1,p}\left(  B_{r_{\Omega}}(x_{0})\right)  .\label{bds}%
\end{equation}

The convergence of $\gamma_{q}$ combined with the boundedness of $u_{q}$
(recall that $\left\Vert u_{q}\right\Vert _{\infty}=1$) and Theorem \ref{Lie}
provides the estimate $\left\Vert u_{q}\right\Vert _{C^{1,\tau}\left(
\overline{\Omega}\right)  }\leq c,$ where the constants $c>0$ and $\tau
\in\left(  0,1\right)  $ are uniform with respect to $q.$ Consequently, up to
another subsequence,
\[
\lim_{q\rightarrow p^{-}}u_{q}=u_{a}\text{\ \ in }C^{1}(\overline{\Omega}),
\]
for some $u_{a}\in W_{0}^{1,p}(\Omega)\cap C^{1}(\overline{\Omega}).$ The
convergence in $C^{1}(\overline{\Omega})$ ensures that we can pass to the
limit (as $q\rightarrow p^{-}$) in (\ref{weak1}) to conclude that $u_{a}$ is a
weak solution to the problem (\ref{Ppw}) and that $\left\Vert u_{a}\right\Vert
_{\infty}=1.$

If $a\geq0$, then it is easy to verify that $\vartheta_{q}$ and $u_{q,\mu_{1}%
}$ form an ordered pair of sub-supersolutions to (\ref{aq}). Hence, Theorem
\ref{sublin} allows us to replace the second inequality in \eqref{bds} with
$\lambda_{a}<\lambda_{1,p}(\Omega)$. The strictness of this inequality is
consequence of \cite[Theorem 2.4]{AH}. Similarly, if $a\leq0,$ then
$\underline{u}_{q}$ and $\vartheta_{q}$ form an ordered pair of
sub-supersolutions to (\ref{aq}), and the first inequaliy in \eqref{bds} can
be replaced with $\lambda_{1,p}(\Omega)\leq\lambda_{a}$.
\end{proof}

\section{Asymptotic behavior as $p\rightarrow\infty\label{secasy}$}

In this section we assume that the coefficient $a$ may depend on $p$ and use
the notation $a_{p}$ to emphasize this\ fact.

We recall that the distance function to the boundary $\partial
\Omega$ is defined as 
\[
d_{\Omega}(x)  =\min_{y\in\partial\Omega}\left\vert
x-y\right\vert \text{ for each }x\in\partial\Omega.
\]
Note that $d_{\Omega}\in C_{0}(\overline{\Omega})\cap W^{1,\infty}(\Omega)$
and $\left\Vert d_{\Omega}\right\Vert _{\infty}=r_{\Omega}.$

In \cite{JLM}, Juutinen, Lindqvist, and Manfredi used the notation
$\Lambda_{\infty}:=r_{\Omega}^{-1},$ which will be kept in this section. This
notation was motivated by the following result they proved: $r_{\Omega}^{-1}$
is the smallest eigenvalue of the $\infty$-Laplacian of the domain $\Omega.$
In the proof of this result, they showed that%
\begin{equation}
\lim_{p\rightarrow\infty}\lambda_{1,p}(\Omega)^{1/p}=\Lambda_{\infty
},\label{JLM1}%
\end{equation}
and deduced a general argument that we state in the next lemma. %A detailed
%proof can be found in \cite[Lemma 2.5]{GE}.

\begin{lemma}
\label{F1}Let $\Omega$ be a bounded domain in $\mathbb{R}^{N}$ and, for each
$p>N$ sufficiently large, let $u_{p}\in W_{0}^{1,p}(\Omega)$ be nonnegative.
If%
\begin{equation}
\limsup_{p\rightarrow\infty}\left\Vert \nabla u_{p}\right\Vert _{p}\leq
C<\infty, \label{gradC}%
\end{equation}
then, up to a subsequence, $u_{p}$ converges uniformly in $\overline{\Omega}$
to a function $u_{\infty}\in C_{0}(\overline{\Omega})\cap W^{1,\infty}%
(\Omega)$ such that
\[
0\leq u_{\infty}(x)\leq\left\Vert \nabla u_{\infty}\right\Vert _{\infty
}d_{\Omega}(x)\leq Cd_{\Omega}(x)\text{ \ for all }x\in\overline{\Omega}.
\]

\end{lemma}

The next lemma can be found in \cite[Theorem 3.2]{EP16}.

\begin{lemma}
\label{F2}Let $\Omega$ be a bounded domain and let
\[
\Lambda_{p}:=\inf\left\{  \left\Vert \nabla u\right\Vert _{p}:u\in W_{0}%
^{1,p}(\Omega)\text{ \ and \ }\left\Vert u\right\Vert _{\infty}=1\right\}  .
\]
Then, $\lim\limits_{p\rightarrow\infty}\Lambda_{p}=\Lambda_{\infty}.$
\end{lemma}

\begin{proposition}
\label{inrad}Let $\lambda_{a_{p}}$ and $u_{a_{p}}\in W_{0}^{1,p}(\Omega)\cap
C^{1}(\overline{\Omega})$ be given by Theorem \ref{teomain}. If
\begin{equation}
\limsup_{p\rightarrow\infty}\left\Vert a_{p}^{+}\right\Vert _{\infty}%
<\Lambda_{\infty}\text{ \ and \ }\lim_{p\rightarrow\infty}\frac{\left\Vert
a_{p}^{-}\right\Vert _{\infty}}{p}=0, \label{hap}%
\end{equation}
then%
\begin{equation}
\lim_{p\rightarrow\infty}\lambda_{a_{p}}^{1/p}=\lim_{p\rightarrow\infty}%
\frac{\left\Vert \nabla u_{a_{p}}\right\Vert _{p}}{\left\Vert u_{a_{p}%
}\right\Vert _{p}}=\Lambda_{\infty}. \label{lim}%
\end{equation}

\end{proposition}

\begin{proof}
We recall that $u_{a_{p}}$ is a positive weak solution to the problem
\begin{equation}
\left\{
\begin{array}
[c]{lll}%
-\Delta_{p}u=\lambda_{a_{p}}\left\vert u\right\vert ^{p-2}u+a_{p}(x)\left\vert
\nabla u\right\vert ^{p-1} & \mathrm{in} & \Omega\\
u=0 & \mathrm{on} & \partial\Omega,
\end{array}
\right.  \label{Pwap}%
\end{equation}
and also that (already using the notation $\Lambda_{\infty}:=r_{\Omega}^{-1}%
$)
\begin{equation}
0<\lambda_{a_{p}}\leq e^{\frac{\left\Vert a_{p}^{-}\right\Vert _{\infty}%
}{\Lambda_{\infty}}}\lambda_{1,p}\left(  B_{r_{\Omega}}(x_{0})\right)
.\label{upb}%
\end{equation}
The inradius of $\Omega$ coincides with the inradius of $B_{r_{\Omega}}%
(x_{0}).$ Thus, it follows from (\ref{JLM1}) that
\[
\lim_{p\rightarrow\infty}\lambda_{1,p}(\Omega)^{1/p}=\lim_{p\rightarrow\infty
}\lambda_{1,p}\left(  B_{r_{\Omega}}(x_{0})\right)  ^{1/p}=\Lambda_{\infty}.
\]
Hence, (\ref{upb}) and the second hypothesis in (\ref{hap}) yield%
\begin{equation}
\limsup_{p\rightarrow\infty}\lambda_{a_{p}}^{1/p}\leq\lim_{p\rightarrow\infty
}e^{\frac{\left\Vert a_{p}^{-}\right\Vert _{\infty}}{p\Lambda_{\infty}}%
}\lambda_{1,p}\left(  B_{r_{\Omega}}(x_{0})\right)  ^{1/p}=\Lambda_{\infty
}.\label{limS}%
\end{equation}

The first hypothesis in (\ref{hap}) implies that%
\[
\limsup_{p\rightarrow\infty}\left\Vert a_{p}^{+}\right\Vert _{\infty}%
\lambda_{1,p}(\Omega)^{-1/p}=\frac{1}{\Lambda_{\infty}}\limsup_{p\rightarrow
\infty}\left\Vert a_{p}^{+}\right\Vert _{\infty}<1.
\]
Consequently, there exists $\theta>0$ such that
\[
\theta\leq\left(  1-\left\Vert a_{p}^{+}\right\Vert _{\infty}\lambda
_{1,p}(\Omega)^{-1/p}\right)  \leq1\text{ for all }p\text{ sufficiently
large.}%
\]

As $u_{a_{p}}$ is a positive weak solution to (\ref{Pwap}) we have that%
\begin{align*}
\left\Vert \nabla u_{a_{p}}\right\Vert _{p}^{p}  &  =\lambda_{a_{p}}\left\Vert
u_{a_{p}}\right\Vert _{p}^{p}+%
%TCIMACRO{\dint \nolimits_{\Omega}}%
%BeginExpansion
{\displaystyle\int\nolimits_{\Omega}}
%EndExpansion
a_{p}(x)\left\vert \nabla u_{a_{p}}\right\vert ^{p-1}u_{a_{p}}\mathrm{d}x\\
&  \leq\lambda_{a_{p}}\left\Vert u_{a_{p}}\right\Vert _{p}^{p}+%
%TCIMACRO{\dint \nolimits_{\Omega}}%
%BeginExpansion
{\displaystyle\int\nolimits_{\Omega}}
%EndExpansion
a_{p}^{+}(x)\left\vert \nabla u_{a_{p}}\right\vert ^{p-1}u_{a_{p}}\mathrm{d}x.
\end{align*}
Hence, H\"{o}lder inequality and the inequality $\left\Vert u_{a_{p}%
}\right\Vert _{p}\leq\left\Vert \nabla u_{a_{p}}\right\Vert _{p}\lambda
_{1,p}(\Omega)^{-1/p}$ yield%
\begin{align*}
\left\Vert \nabla u_{a_{p}}\right\Vert _{p}^{p}  &  \leq\lambda_{a_{p}%
}\left\Vert u_{a_{p}}\right\Vert _{p}^{p}+\left\Vert a_{p}^{+}\right\Vert
_{\infty}\left\Vert u_{a_{p}}\right\Vert _{p}\left\Vert \nabla u_{a_{p}%
}\right\Vert _{p}^{p-1}\\
&  \leq\lambda_{a_{p}}\left\Vert u_{a_{p}}\right\Vert _{p}^{p}+\left\Vert
a_{p}^{+}\right\Vert _{\infty}\lambda_{1,p}(\Omega)^{-1/p}\left\Vert \nabla
u_{a_{p}}\right\Vert _{p}^{p},
\end{align*}
so that%
\[
\frac{\left\Vert \nabla u_{a_{p}}\right\Vert _{p}^{p}}{\left\Vert u_{a_{p}%
}\right\Vert _{p}^{p}}\left(  1-\left\Vert a_{p}^{+}\right\Vert _{\infty
}\lambda_{1,p}(\Omega)^{-1/p}\right)  \leq\lambda_{a_{p}}%
\]
It follows that
\begin{equation}
\lambda_{1,p}(\Omega)^{1/p}\leq\frac{\left\Vert \nabla u_{a_{p}}\right\Vert
_{p}}{\left\Vert u_{a_{p}}\right\Vert _{p}}\leq\frac{\lambda_{a_{p}}^{1/p}%
}{\theta^{1/p}}, \label{asy1}%
\end{equation}
for $p$ large enough. Therefore, taking into account that $\lim_{p\rightarrow
\infty}\theta^{1/p}=1$ we obtain (\ref{lim}) from (\ref{asy1}) and (\ref{limS}).
\end{proof}

Owing to the homogeneity of (\ref{Pwap}), we can assume that $u_{a_{p}}$ is
$L^{p}$-normalized, that is, $\left\Vert u_{a_{p}}\right\Vert _{p}=1.$

\begin{corollary}
In addition to the hypotheses of Proposition \ref{inrad}, assume that
$\left\Vert u_{a_{p}}\right\Vert _{p}=1.$ Then,
\begin{equation}
\lim_{p\rightarrow\infty}\left\Vert u_{a_{p}}\right\Vert _{\infty}=1,
\label{asy2}%
\end{equation}
and, up to a subsequence, $u_{a_{p}}$ converges uniformly, as $p\rightarrow
\infty,$ to a function $u_{\infty}\in C_{0}(\overline{\Omega})\cap
W^{1,\infty}(\Omega)$ satisfying the properties:

\begin{enumerate}
\item $0\leq u_{\infty}(x)\leq\Lambda_{\infty}d_{\Omega}(x)$ for all $x\in$
$\overline{\Omega},$

\item $\left\Vert u_{\infty}\right\Vert _{\infty}=1$ \ and \ $\Lambda_{\infty
}=\left\Vert \nabla u_{\infty}\right\Vert _{\infty},$

\item $\left\{  x\in\Omega:u_{\infty}(x)=1\right\}  \subseteq\left\{
x\in\Omega:d_{\Omega}(x)=\left\Vert d_{\Omega}\right\Vert _{\infty}\right\}
.$
\end{enumerate}
\end{corollary}

\begin{proof}
Let us denote the $N$-Lebesgue measure of $\Omega$ by $\left\vert
\Omega\right\vert .$ We have that%
\[
1=\left\Vert u_{a_{p}}\right\Vert _{p}\leq\left\Vert u_{a_{p}}\right\Vert
_{\infty}\left\vert \Omega\right\vert ^{1/p}\text{ \ and \ }\left\Vert
u_{a_{p}}\right\Vert _{\infty}\leq\Lambda_{p}^{-1}\left\Vert \nabla u_{a_{p}%
}\right\Vert _{p}.
\]
Hence, letting $p\rightarrow\infty$ in the inequalities%
\[
\left\vert \Omega\right\vert ^{-1/p}\leq\left\Vert u_{a_{p}}\right\Vert
_{\infty}\leq\Lambda_{p}^{-1}\left\Vert \nabla u_{a_{p}}\right\Vert _{p}%
\]
we obtain (\ref{asy2}) from Lemma \ref{F2} and (\ref{lim}).

Now, let us prove items 1,2, and 3 all at once. According to Lemma \ref{F1},
up to a subsequence, $u_{a_{p}}$ converges uniformly in $\overline{\Omega}$ to
a function $u_{\infty}\in C_{0}(\overline{\Omega})\cap W^{1,\infty}(\Omega)$
such that
\[
0\leq u_{\infty}(x)\leq\left\Vert \nabla u_{\infty}\right\Vert _{\infty
}d_{\Omega}(x)\leq\Lambda_{\infty}d_{\Omega}(x)\text{ \ for all }x\in
\overline{\Omega}.
\]

The uniform convergence implies that
\[
\left\Vert u_{\infty}\right\Vert _{\infty}=\lim_{p\rightarrow\infty}\left\Vert
u_{a_{p}}\right\Vert _{\infty}=1.
\]
Thus, if $x_{0}\in\Omega$ is such that $u_{\infty}(x_{0})=1,$ then%
\[
1\leq\left\Vert \nabla u_{\infty}\right\Vert _{\infty}d_{\Omega}(x_{0}%
)\leq\Lambda_{\infty}d_{\Omega}(x_{0})\leq1.
\]
These inequalities imply that $d_{\Omega}(x_{0})=\Lambda_{\infty}%
^{-1}=\left\Vert \nabla u_{\infty}\right\Vert _{\infty}^{-1}.$
\end{proof}

\appendix 
\section [Appendix A] {Appendix: uniqueness of the trivial solution when $\lambda \leq 0$} \label{null}

\begin{proposition}
If $\lambda \leq 0$, then the only nonnegative weak solution to \eqref{Pw} is the trivial one. 
\end{proposition}

\begin{proof}
For each $t\in\lbrack0,M),$ where $M:=\operatorname*{ess}\sup u,$ define
$A_{t}:=\left\{  x\in\Omega:u>t\right\}  $ and $(u-t)^{+}:=\max\left\{
u-t,0\right\}  \in W_{0}^{1,p}(\Omega).$ Since $\nabla(u-t)=\nabla u$ in
$A_{t},$ the well-known Poincar\'{e}-Sobolev inequality yields%
\begin{equation}%
%TCIMACRO{\dint \nolimits_{A_{t}}}%
%BeginExpansion
{\displaystyle\int\nolimits_{A_{t}}}
%EndExpansion
(u-t)^{p}\mathrm{d}x\leq C_{N,p}\left\vert A_{t}\right\vert ^{p/N}%
%TCIMACRO{\dint \nolimits_{A_{t}}}%
%BeginExpansion
{\displaystyle\int\nolimits_{A_{t}}}
%EndExpansion
\left\vert \nabla u\right\vert ^{p}\mathrm{d}x,\label{PS}%
\end{equation}
where $C_{N,p}$ is a positive constant that depends only on $p$ and $N.$

Moreover, as%
\[
0\leq%
%TCIMACRO{\dint \nolimits_{\Omega}}%
%BeginExpansion
{\displaystyle\int\nolimits_{\Omega}}
%EndExpansion
(u-t)^{+}\mathrm{d}x=%
%TCIMACRO{\dint \nolimits_{A_{t}}}%
%BeginExpansion
{\displaystyle\int\nolimits_{A_{t}}}
%EndExpansion
(u-t)\mathrm{d}x=%
%TCIMACRO{\dint \nolimits_{A_{t}}}%
%BeginExpansion
{\displaystyle\int\nolimits_{A_{t}}}
%EndExpansion
u\mathrm{d}x-t\left\vert A_{t}\right\vert \leq\left\Vert u\right\Vert
_{1}-t\left\vert A_{t}\right\vert ,
\]
we have that
\begin{equation}
t\left\vert A_{t}\right\vert \leq\left\Vert u\right\Vert _{1}.\label{At0}%
\end{equation}

Using $(u-t)^{+}$ in the weak formulation of \eqref{Pw} we obtain
\begin{align*}%
%TCIMACRO{\dint \nolimits_{A_{t}}}%
%BeginExpansion
{\displaystyle\int\nolimits_{A_{t}}}
%EndExpansion
\left\vert \nabla u\right\vert ^{p}\mathrm{d}x  & =\lambda%
%TCIMACRO{\dint \nolimits_{A_{t}}}%
%BeginExpansion
{\displaystyle\int\nolimits_{A_{t}}}
%EndExpansion
u^{p-1}(u-t)\mathrm{d}x+%
%TCIMACRO{\dint \nolimits_{A_{t}}}%
%BeginExpansion
{\displaystyle\int\nolimits_{A_{t}}}
%EndExpansion
a(x)\left\vert \nabla u\right\vert ^{p-1}(u-t)\mathrm{d}x\\
& \leq\left\Vert a^{+}\right\Vert _{\infty}%
%TCIMACRO{\dint \nolimits_{A_{t}}}%
%BeginExpansion
{\displaystyle\int\nolimits_{A_{t}}}
%EndExpansion
\left\vert \nabla u\right\vert ^{p-1}(u-t)\mathrm{d}x.
\end{align*}
Hence, H\"{o}lder inequality and
(\ref{PS}) yield%
\begin{align*}%
%TCIMACRO{\dint \nolimits_{A_{t}}}%
%BeginExpansion
{\displaystyle\int\nolimits_{A_{t}}}
%EndExpansion
\left\vert \nabla u\right\vert ^{p}\mathrm{d}x &  \leq\left\Vert
a^{+}\right\Vert _{\infty}\left(
%TCIMACRO{\dint \nolimits_{A_{t}}}%
%BeginExpansion
{\displaystyle\int\nolimits_{A_{t}}}
%EndExpansion
\left\vert \nabla u\right\vert ^{p}\mathrm{d}x\right)  ^{\frac{p-1}{p}}\left(
%
%TCIMACRO{\dint \nolimits_{A_{t}}}%
%BeginExpansion
{\displaystyle\int\nolimits_{A_{t}}}
%EndExpansion
(u-t)^{p}\mathrm{d}x\right)  ^{\frac{1}{p}}\\
&  \leq\left\Vert a^{+}\right\Vert _{\infty}\left(
%TCIMACRO{\dint \nolimits_{A_{t}}}%
%BeginExpansion
{\displaystyle\int\nolimits_{A_{t}}}
%EndExpansion
\left\vert \nabla u\right\vert ^{p}\mathrm{d}x\right)  ^{\frac{p-1}{p}}%
C_{N,p}^{1/p}\left\vert A_{t}\right\vert ^{1/N}\left(
%TCIMACRO{\dint \nolimits_{A_{t}}}%
%BeginExpansion
{\displaystyle\int\nolimits_{A_{t}}}
%EndExpansion
\left\vert \nabla u\right\vert ^{p}\mathrm{d}x\right)  ^{\frac{1}{p}},
\end{align*}
so that%
\begin{equation}%
%TCIMACRO{\dint \nolimits_{A_{t}}}%
%BeginExpansion
{\displaystyle\int\nolimits_{A_{t}}}
%EndExpansion
\left\vert \nabla u\right\vert ^{p}\mathrm{d}x\leq C_{N,p}^{1/p}\left\Vert
a^{+}\right\Vert _{\infty}\left\vert A_{t}\right\vert ^{1/N}%
%TCIMACRO{\dint \nolimits_{A_{t}}}%
%BeginExpansion
{\displaystyle\int\nolimits_{A_{t}}}
%EndExpansion
\left\vert \nabla u\right\vert ^{p}\mathrm{d}x.\label{simp}%
\end{equation}

Now, let us suppose, by contradiction, that $u\not \equiv 0.$ It follows that
$M>0$ and%
\[%
%TCIMACRO{\dint \nolimits_{A_{t}}}%
%BeginExpansion
{\displaystyle\int\nolimits_{A_{t}}}
%EndExpansion
(u-t)^{p}\mathrm{d}x>0\text{ for all }t\in[0,M).
\]
Hence, according to (\ref{PS}) we can divide (\ref{simp}) by $%
%TCIMACRO{\dint \nolimits_{A_{t}}}%
%BeginExpansion
{\displaystyle\int\nolimits_{A_{t}}}
%EndExpansion
\left\vert \nabla u\right\vert ^{p}\mathrm{d}x$ and use (\ref{At0}) to achieve
the estimate%
\[
t\leq C_{N,p}^{N/p}\left\Vert a^{+}\right\Vert _{\infty}^{N}\left\Vert
u\right\Vert _{1},\text{ for all }t\in(0,M).
\]
This means that $u\in L^{\infty}(\Omega)$ and $M=\left\Vert u\right\Vert
_{\infty}\leq C_{N,p}^{N/p}\left\Vert a^{+}\right\Vert _{\infty}^{N}\left\Vert
u\right\Vert _{1}.$

By Theorem \ref{Lie}, $u\in C^{1}(\overline{\Omega}).$ It follows that
$u(x_{0})=M$ for some $x_{0}\in\Omega.$ Let us define $v\in W_{0}^{1,p}%
(\Omega)\cap C^{1}(\overline{\Omega})$ by $v(x):=M-u(x)$ for all
$x\in\overline{\Omega}.$ This function is nonnegative and satisfies
\[
\Delta_{p}v-a(x)\left\vert \nabla v\right\vert ^{p-1}\leq0\text{ \ in }\Omega,
\]
in the distribution sense. As $v(x_{0})=0,$ we can apply the Strong Maximum
Principle by Pucci, Serrin and Zou \cite[Theorem 1']{PSZ}, with $k=\left\Vert
a^{+}\right\Vert _{\infty}$ and $f\equiv0$ (see also \cite[Theorem 5.3.1]{PS}), to conclude that $v\equiv0.$ It
follows that $u\equiv M$. However, the continuity of $u$ up to
the boundary contradicts the positivity of $M$. Therefore, $u\equiv0.$
\end{proof}

%\section*{Acknowledgements}

%The author was partially supported by Fapemig/Brazil (RED-00133-21) and
%CNPq/Brazil (307866/2025-4).

\section*{Declarations}

\noindent {\bf Acknowledgements.} The author was partially supported by Fapemig/Brazil (RED-00133-21) and
CNPq/Brazil (307866/2025-4).
\newline

\noindent {\bf Competing Interests.} The author declares that there is no conflict of interests.
\newline

\noindent {\bf Data Availability.}  Data availability is not applicable to this article as no datasets were generated or analyzed during the current study.


\begin{thebibliography}{99}                                                                                               %


\bibitem {AP}B. Abdellaoui and I. Peral: Existence and nonexistence results
for quasilinear elliptic equations involving the \textit{p}-Laplacian with a
critical potential, Ann. Mat. Pura Appl. \textbf{182} (2003) 247--270.

\bibitem {AH}W. Allegreto and Y.X. Huang: A Picone's identity for the
$p$-Laplacian and applications, Nonlinear Anal. \textbf{32} (1998) 819--830.

\bibitem {ACG}A. Anane, O. Chakrone and J. P. Gossez: Spectre d'ordre
sup\'{e}rieur et probl\`{e}mes aux limites quasi-lin\'{e}aires, Bollettino
dell'Unione Matematica Italiana, \textbf{4B} (2001) 483-519.

\bibitem {BJEM}R.J. Biezuner, J. Brown, G. Ercole and E.M. Martins: Computing the
first eigenpair of the p-Laplacian via inverse iteration of Sublinear
Supersolutions, J. Sci. Comput. \textbf{52} (2012) 180--201.

\bibitem {BMP}L. Boccardo, F. Murat and J.-P. Puel: R\'{e}sultats
d'existence pour certains probl\`{e}mes elliptiques quasilin\'{e}aires, Ann.
Scuola Norm. Sup. Pisa Cl. Sci. \textbf{11} (1984) 213--235.

\bibitem {CL}M. Cuesta Leon: Existence results for quasilinear problems via
ordered sub and supersolutions, Ann. Fac. Sci. Toulouse \textbf{6} (1997) 591--608.

\bibitem {DSS} Y. de la Calle, J.C. Sabina de Lis and S. Segura de Leon: The convective eigenvalues of the one-dimensional $p$-Laplacian as $p\to 1$,  J. Math. Anal. Appl. \textbf{484} (2020) 123738.

%\bibitem {GE}G. Ercole: Asymptotics for Sobolev extremals: the hyperdiffusive
%case, Ann. Mat. Pura Appl. (to appear) https://doi.org/10.1007/s10231-026-01717-w.

\bibitem {EP16}G. Ercole and G.A. Pereira: Asymptotics for the best Sobolev
constants and their extremal functions, Math. Nachr. \textbf{289} (2016) 1433--1449.

\bibitem {GST}J. Garc\'{\i}a-Mel\'{\i}an, J.C. Sabina de Lis and P.
Tak\'{a}\v{c}: Dirichlet problems for the $p$-Laplacian with a convection
term, Rev. Mat. Complut. \textbf{30} (2017) 313--334.

\bibitem{GSS} J. Garc\'{\i}a-Mel\'{\i}an, J.C. Sabina de Lis and M.  Sanabria-Garc\'{\i}a:  Eigenvalue analysis for the p-Laplacian under convective perturbation. J. Comput. Appl. Math. 110, 73–91 (1999)

\bibitem {JLM}P. Juutinen, P. Lindqvist and J. Manfredi: The $\infty$-eigenvalue
problem, Arch. Ration. Mech. Anal. \textbf{148} (1999) 89--105.

\bibitem {Lieb}G.M. Lieberman: Boundary regularity for solutions of degenerate
elliptic equations, Nonlinear Anal. \textbf{12} (1988) 1203--1219.

\bibitem {PS}P. Pucci and J. Serrin. The maximum principle. Basel: Birkhäuser Basel, 2007.

\bibitem {PSZ} P. Pucci, J. Serrin and H. Zou: A strong maximum principle and a compact support principle for singular elliptic inequalities,  J. Math. Pures Appl., \textbf{78} (1999) 769-–789.

\bibitem {SP}H. Shouchuan and N.S. Papageorgiou: Positive solutions for nonlinear
Dirichlet problems with convection, Appl. Math. Optim. \textbf{82} (2020) 451--470.
\end{thebibliography}
\end{document}